\documentstyle[11pt]{article}

\begin{document}

\noindent

\begin{center}

  {\LARGE Stringy  orbifolds}\footnote{ partially
supported by the National Science Foundation and Hong Kong RGC grant}

  \end{center}

  \noindent

  \begin{center}

     {\large  Yongbin Ruan}\\[5pt]

      Department of Mathematics, University of Wisconsin-Madison\\

        and Hong Kong University of Science and Technology\\[5pt]

        Email:  ruan@math.wisc.edu and ruan@ust.hk\\[5pt]

              \end{center}

              \def \x{{\bf x}}

              \def \J{{\cal J}}

              \def \M{{\cal M}}

              \def \A{{\cal A}}

              \def \B{{\cal B}}

              \def \C{{\bf C}}

              \def \Z{{\bf Z}}

              \def \R{{\bf R}}
          \def \RR{{\cal R}}
          \def \UU{{\cal U}}
              \def \P{{\bf P}}
              \def \LL{{\cal L}}
                \def \GL{{\cal GL}}
                \def \GG{{\cal G}}
              \def \I{{\bf I}}
                \def \HH{{\cal H}}
                \def \Xx{{\sf X}}
                \def \FF{{\cal F}}
              \def \N{{\cal N}}

              \def \T{{\cal T}}

              \def \Q{{\bf Q}}

              \def \D{{\cal D}}

              \def \H{{\cal H}}

              \def \S{{\cal S}}

              \def \e{{\bf E}}

              \def \CP{{\bf CP}}

              \def \U{{\cal U}}

              \def \E{{\cal E}}

              \def \F{{\bf C}}

              \def \L{{\cal L}}

              \def \K{{\cal K}}

              \def \G{{\cal G}}

              \def \z{{\bf z}}

              \def \m{{\bf m}}

              \def \n{{\bf n}}

              \def \V{{\cal V}}

              \def \W{{\cal W}}

 \def \g{{\bf g}}

              \def \h{{\bf h}}

                \def \O{{\cal O}}

\tableofcontents

\section{Introduction}
    During the last two years, there was an explosive growth of activities to study the "stringy" properties of
orbifolds, motivated by the orbifold stringy theory from physics
in the middle eighty. A year ago, this author wrote an article to
survey some of these initial developments. In the May of 2000, a
very successful conference on orbifold was held in Madison where
more than 40 mathematician and physicists came together to talk
about orbifolds. It was clear that the vast content of this new
subject of mathematics is much beyond the author's expectation.
Hence, his initial survey was completely outdated. A new
survey is needed to expose the new connections and new results the
author leaned on and after conference to general audience. This is
the main motivation and purpose of this survey article.
Originally, the author planned to give a comprehensive survey to
cover the most important results and bring the reader the frontier
of research. Because of the limited time, the author has to scale
down his plan to cover the only topics he is most familiar with.
He  will list other topics in the last sections.

    The article consists of several relatively independent topics. Let me briefly sketch its content. On the Madison conference,
    the author leant the
    interpretation of the orbifold in terms of groupoid by the work of Moerdijk and his collaborators \cite{CM}, \cite{MP1},\cite{MP2}.
    This is an important
    approach to orbifolds. This approach will be presented in the section 2. The groupoid approach  is a categorical and more abstract approach.
    Although it is less concrete, but it does looks cleaner.
Besides, Lupercio-Uribe's twisted K-theory seems to be best
described in this framework. Orbifold and groupoid provided two
different points of view. The author feels that it is important to
have diverse points of view on this subject. The work of Lupercio
and
    Uribe \cite{LU1}, \cite{LU2} on the twisted K-theory suggests that many of results in orbifold should be pushed to more general orbispace
    where the local model is the quotient of a smooth manifold by a Lie group. This kind of space is called Artin stack and also
    admitted a groupoid description. However, in many way, Satake's original approach is more concrete and easier to work with
    in differential geometry. In the section 2, we present the version of orbispace and its morphism by W. Chen \cite{C} which follows
    Satake's original approach.

    In the section 3, we will focus on the algebraic geometrical application of stringy orbifold. Here, we focus on the
computation of the cohomological ring structure of the crepant resolution
of orbifolds, a topic attracting many algebraic geometors in the
past. The main example we use is the Hilbert scheme of points of
the algebraic surfaces. In this section, we will present the beautiful
calculation of the orbifold cohomology of the symmetry product of the complex
manifolds by Vafa-Witten (vector space) \cite{VW},
Fantechi-Gottsche-Uribe (ring structure) \cite{FG2}, \cite{U}.
Together with Lehn-Sorger's computation of the cohomology of Hilbert
scheme of points of $K3$ and $T^4$ \cite{LS2}. They verified a
conjecture of the author for the hyperkahler resulotion of orbifolds.
In the end, we will present the author's conjecture for general crepant
resolution and verify it for a few examples.

    An important aspect of stringy orbifolds is the (twisted)-orbifold K-theory \cite{AR}. Soon after the work of Adem-Ruan,
    the twisted orbifold K-theory and the twisted K-theory on smooth
    manifolds by \cite{WIT2} were unified by Lupercio and
Uribe \cite{LU1} as the special cases of the twisted K-theory on
groupoid by gerbes. We present their construction in the section
4. In the light of recent computation of equivariant twisted
K-theory of a Lie group acting on itself by the conjugation by
Freed-Hopkins-Teleman \cite{FR}, many of current theories should be
pushed to the more general orbispace where an isotropy subgroup is
not necessarily finite. There are many questions here. It seems to
be a very exciting research direction in next few years. In this
section, we will also present an interesting construction of Wang
to  twist the K-theory using spin represenation (fermion) \cite{W}. As
Wang showed for the symmetric product, its Euler characteristic
formula is much nicer than that of the twisted K-theory of Adem-Ruan.
This suggests that more should be done regarding to Wang's
twisting.

    One of attractions of the stringy orbifold is its unique orbifold feature, which does not exist on smooth manifolds.
    Orbifold cohomology is such an example. So far, the most of actions are on classical theories such as cohomology and K-theory.
    Now, we add an unique orbifold feature to its quantum theory. This is the integration of the theory of spin curves by Jarvis and
    others \cite{J1} \cite{J2} \cite{JKV} \cite{JKV2}
into orbifolds. The theory of spin curves was motivated by the
2D-gravity coupled with the matter and has been around for ten
years \cite{WIT1}.
 In my view, it is
a piece of mathematics before its time. There is no doubt that its
natural home is in the stringy orbifold \cite{AJ}. Furthermore, it
also provides a missing link in the orbifold quantum cohomology.
Recall that physically GW-invariants are correlation functions of
2D-gravities coupled with topological sigma model. The correct
analogous of the topological sigma model in orbifold is not an ordinary
map but a good map/groupoid morphism. It is amazing that the correct
analogous of 2D-gravity in the orbifold is not stable curves either.
Instead, it is a spin curve or 2D-gravity coupled with the matter.
From this point of view, it gives a satisfactory construction of the
(spin) orbifold Gromov-Witten invariants over orbifolds. It is not
clear at this moment for its implication in other topics of
interest such as mirror symmetry. For example, what is the B-model
correspondence of the spin orbifold quantum cohomology? In this
section, we will sketch Jarvis-Ruan's construction \cite{JR} of
the spin GW-invariants.

    In last sections, we make some general remarks about other
    important topics in orbifolds.

    This author would like to thank all the participants to the Workshop on Mathematical Aspect of Orbifold Stringy Theory and their
wonderful talks where the original motivation of this paper was conceived. In particular, he would like to thank his co-organizers
A. Adem and J. Morava for the wonderful job they have done so that the current organizer doesn't have to organize at all. He also would like to thank
E. Lupercio for his unyielding effort to educate him for the theory of groupoid and stack.

\section{Foundation}
    An
important aspect of the foundation of orbifold came to my attention
during the Madison conference is the interpretation of the orbifold as
a groupoid through the earlier work of Moerdijk and his collaborators
\cite{CM}\cite{MP1}\cite{MP2}. Compared to Satake's approach,
the groupoid is a categorical notion and more abstract. For my
personal taste, I will probabily prefer Satake's ad hoc approach
because it is more concrete and easy to work with in the differential geometry.
But the more abstract approach of groupoid does have its
advantage. First of all, it connects to the theory of the stack in the
algebraic geometry. Secondly, many important constructions such as
the twisted K-theory by gerbes \cite{LU1} is best described using this
framework. The author feels that it is useful to keep this dual
approache of the orbifold. Therefore, we will present the groupoid
approach of the orbifold in this section.

Compared to groupoid approach, Satake's original approach looks
like ad hoc. But it has the advantages to be more concrete and easier to
work with in differential geometry. Even though its definition due
to Satake was known since 1950's, somehow it was very much under
developed.
  Basically, it was viewed as a generalized smooth structure and has no distinct character of its own. When the author and
W. Chen developed the Gromov-Witten theory over orbifolds, there
was very little references we could use and we had to start from
scratch. We discovered that Satake's notion of an orbifold map is not good
for pulling back bundles and we cooked out a new notion called {\em
good map}. Recently, it was showed to be equivalent to the
morphism of the groupoid (see \cite{LU1}).  However, it seems to be
easier to do the computation for good maps. For example, in the orbifold
quantum cohomology, one needs to compute the number of good
maps/groupoid morphisms from the orbifold Riemann surface to a
symplectic orbifold which realizes the same map on the underline
topological space. For good maps, we have a concrete
classification. It seems to be harder to do so for groupoid
morphism.

Recently, Lupercio-Uribe unified the twisted orbifold K-theory of
Adem-Ruan with the twisted K-theory on smooth manifold \cite{LU1}.
Moreover, they showed that their construction works for more
general space of the Artin stack. It strongly suggests that the theory
of stringy orbifolds can and should be push to this more general
category. Groupoid approach automatically includes this case. It
should be useful to extend Satake's approach to this case as well.
This was done by W. Chen \cite{C}, which we will present it here.

\subsection{Orbifold and Groupoid}

 By the work of Moerdijk and his collaborators, an orbifold has a nice interpretation as a groupoid.
We shall sketch this connection in this section. A groupoid can be
thought of as a generalization of a group, a manifold and an
equivalence relation. First as an equivalence relation, a groupoid
has a set of relations $\RR$ that we will think of as arrows.
These elements arrows relate  is a set $\UU$. Given an arrow
$\stackrel{r}{\rightarrow} \in \RR$ it has a source
$x=s(\stackrel{r}{\rightarrow})\in\UU$ and a target
$y=t(\stackrel{r}{\rightarrow})\in\UU$. Then we say that
$x\stackrel{r}{\rightarrow} y$, namely $x$ is related to $y$. We
want to have an equivalence relation, for example we want
transitivity and then we will need a way to compose arrows
$x\stackrel{r}{\rightarrow} y\stackrel{r}{\rightarrow}z$. We also
require $\RR$ and $\UU$ to be more than mere sets. Sometimes we
want them to be locally Hausdorff, paracompact, locally compact
topological spaces, smooth manifolds.

\vskip 0.1in \noindent {\bf Definitions 2.1: }{\it A
\emph{groupoid} is a pair of objects in a category $\RR, \UU$ and
morphisms
$$s,t : \RR \stackrel{\rightarrow}{\rightarrow} \UU$$ called respectively source and target ,
provided with an identity $$e : \UU \rightarrow \RR $$ a
multiplication
$$m:\RR_{t}\times_s \RR \rightarrow \RR$$ and an inverse
$$i : \RR \rightarrow \RR$$ satisfying the following properties:

\begin{enumerate}
\item The identity inverts both $s$ and $t$:
$$ \begin{array}{ccccccccccccc}
 \UU& &\stackrel{e}{\rightarrow} & & \RR& & & & \UU& &\stackrel{e}{\rightarrow} & & \RR\\
 & & \stackrel{id_\UU}{\searrow} & &\downarrow s & &and & & & &\stackrel{id_\UU}{\searrow} & &\downarrow t\\
 & & & &\UU & & & & & & & &\UU
\end{array}$$

\item Multiplication is compatible with both $s$ and $t$:

$$ \begin{array}{ccccccccccccc}
 \RR_{t}\times_s \RR & &\stackrel{m}{\rightarrow} & &\RR & & & &\RR_{t}\times_s \RR  & &\stackrel{m}{\rightarrow} & & \RR \\
\downarrow \pi_1 & & & &\downarrow s & &and & &\downarrow \pi_2 & & & & \downarrow t\\
\RR & &\stackrel{s}{\rightarrow} & &\UU & & & &\RR &
&\stackrel{t}{\rightarrow} & & \UU
\end{array}$$

\item Associativity:
$$ \begin{array}{ccccc}
\RR_{t} \times_s \RR_{t} \times_s \RR  & &\stackrel{id_\RR \otimes m}{\rightarrow} & &\RR_{t} \times_s \RR  \\
 \downarrow m \otimes id_\RR& & & & \downarrow m \\
\RR_{t} \times_s \RR & & \stackrel{m}{\rightarrow}& &\RR
\end{array}$$

\item Unit condition:
$$ \begin{array}{ccccccccccccc}
 \RR& &\stackrel{(e \circ s,id_\RR)}{\rightarrow} & & \RR_{t} \times_s \RR& & & & \RR& &\stackrel{(id_\RR,e \circ t)}{\rightarrow}
  & & \RR_{t}\times_s \RR\\
 & & \stackrel{id_\RR}{\searrow} & &\downarrow m & &and & & & &\stackrel{id_\RR}{\searrow} & &\downarrow m\\
 & & & &\RR & & & & & & & &\RR
\end{array}$$

\item Inverse:
$$
\begin{array}{ccc}
i \circ i & = & id_\RR \\
s \circ i & = & t \\
t \circ i & = & s
\end{array}
$$
with

$$ \begin{array}{ccccccccccccc}
\RR & &\stackrel{(i, id_\RR)}{\rightarrow} & &\RR_{t}\times_s \RR & & & & \RR  & &\stackrel{(i,id_\RR)}{\rightarrow} & & \RR_{t}\times_s \RR \\
\downarrow s & & & &\downarrow m & &and & &\downarrow t & & & & \downarrow m\\
\UU & &\stackrel{e}{\rightarrow} & &\RR & & & &\UU &
&\stackrel{e}{\rightarrow} & & \RR
\end{array}$$

\end{enumerate}}

We denote the groupoid by  $\RR
\stackrel{\rightarrow}{\rightarrow} \UU : = (\RR, \UU,
s,t,e,m,i)$, and the groupoid is called \emph{\'etale } if the
base category is that of locally Hausdorff, paracompact, locally
compact topological spaces and the maps $s,t: \RR \rightarrow U$
are local homeomorphisms (diffeomorphisms). We will say that a
groupoid is proper is $s$ and $t$ are proper  maps. From now on we
will assume that our groupoids are differentiable, \'etale and
proper.

\vskip 0.1in \noindent {\bf Example 2.2: }{\it For $M$ a manifold
and $\{U_\alpha\}$ and open cover, let
$$\UU = \bigsqcup_\alpha U_\alpha \ \ \  \ \RR= \bigsqcup_{(\alpha, \beta)} U_\alpha \cap U_\beta \ \ (\alpha,\beta) \neq (\beta,\alpha)$$
$$s|_{U_{\alpha \beta}}: U_{\alpha \beta} \rightarrow U_\alpha, \ t|_{U_{\alpha \beta}}: U_{\alpha \beta} \rightarrow U_\beta \ \
e|_{U_{\alpha}}: U_\alpha \rightarrow U_\alpha$$
 $$i|_{U_{\alpha \beta}}: U_{\alpha \beta} \rightarrow U_{\beta \alpha} \ \& \
m|_{U_{\alpha \beta \gamma}}: U_{\alpha \beta \gamma} \rightarrow
U_{\alpha \gamma}$$ the natural maps. Note that in this example
$\RR_{t} \times_s \RR$ coincides with the subset of $\RR_{t}
\times_s \RR$ of pairs $(u,v)$ so that $t(u)=s(v)$, namely the
disjoint union of all possible triple intersections
$U_{\alpha\beta\gamma}$ of open sets in the open cover
$\{U_\alpha\}$.}

\vskip 0.1in \noindent {\bf Example 2.3: }{\it Let $G$ be a group
and $U$ a set provided with a left $G$ action
$$G \times U \rightarrow U$$
$$ (g,u) \rightarrow gu$$
we put $\UU=U$ and $\RR = G \times U$ with $s(g,u)=u$ and
$t(g,u)=gu$. The domain of $m$ is the same as $G \times G \times
U$ where $m(g,h,u)=(gh,u)$, $i(g,u)=(g^{-1},gu)$ and
$e(u)=(id_G,u)$.  } \vskip 0.1in

We will write $G \times U \stackrel{\rightarrow}{\rightarrow} U$
(or sometimes $X=[U/G]$,) to denote this groupoid. When $U=\star$
is a point, we denote the groupoid by $\bar{G}$.

\vskip 0.1in \noindent {\bf Definition 2.4: }{\it A morphism of
grupoids $(\Psi, \psi): (\RR' \stackrel{\rightarrow}{\rightarrow}
\UU') \rightarrow (\RR \stackrel{\rightarrow}{\rightarrow} \UU)$
are the following commutative diagrams:

$$\left. \begin{array}{ccccc}
 \RR' & & \stackrel{\Psi}{\rightarrow} & & \RR \\
 s' \downarrow\downarrow t' & & & & s \downarrow \downarrow t \\
\UU' & & \stackrel{\psi}{\rightarrow} & & \UU
\end{array} \right.
 \ \ \ \ \ \ \  \ \
\left. \begin{array}{ccccc}
 \RR' & & \stackrel{\Psi}{\rightarrow} & & \RR \\
 e' \uparrow & & & &  e \uparrow \\
\UU' & & \stackrel{\psi}{\rightarrow} & & \UU
\end{array} \right.
$$
$$
\left. \begin{array}{ccccc}
 \RR'_{t'}\times_{s'} \RR' & & \stackrel{\Psi}{\rightarrow} & & \RR_{t} \times_s \RR \\
 m' \downarrow & & & & m \downarrow \\
\RR' & & \stackrel{\psi}{\rightarrow} & & \RR
\end{array} \right.
 \ \ \ \ \ \ \
\left. \begin{array}{ccccc}
 \RR' & & \stackrel{\Psi}{\rightarrow} & & \RR \\
 i' \downarrow & & & & i \downarrow \\
\RR' & & \stackrel{\psi}{\rightarrow} & & \RR
\end{array} \right.$$

} Now we need to say when two groupoids are ``equivalent''

\vskip 0.1in \noindent {\bf Definition 2.5: }{\it A morphism of
groupoids $(\Psi,\psi)$ is called \emph{a Morita morphism}
whenever:
\begin{itemize}
\item The map $s \circ \pi_2 : \UU'_{\psi}\times_{t} \RR \rightarrow \UU$ is an \'etale surjection,
\item The following square is a fibered product
$$\begin{array}{ccccc}
 \RR' & & \stackrel{\Psi}{\rightarrow} & & \RR \\
 (s',t')\downarrow \downarrow & & & & (s,t) \downarrow \downarrow\\
\UU'\times \UU' & & \stackrel{\psi \times \psi}{\rightarrow} & &
\UU \times \UU
\end{array}$$
\end{itemize}}
\vskip 0.1in

Two groupoids $\RR_1 \stackrel{\rightarrow}{\rightarrow} \UU_1$,
$\RR_2 \stackrel{\rightarrow}{\rightarrow} \UU_2$ are called
\emph{Morita equivalent} if there are Morita morphisms $(\Psi_i,
\psi_i): \RR' \stackrel{\rightarrow}{\rightarrow} \UU' \rightarrow
\RR_i \stackrel{\rightarrow}{\rightarrow} \UU_i$ for $i = 1,2$.
This is an equivalence relation and in general we will consider
the category of \'etale groupoids obtained by formally inverting
the Morita equivalences (see \cite{Moerdijk} for details).

 Let $X$ be an orbifold and
$\{(V_p,G_p,\pi_p)\}_{p \in X}$ its orbifold structure, the
groupoid $\RR \stackrel{\rightarrow}{\rightarrow} \UU$ associated
to $X$ will be defined as follows: $\UU := \bigsqcup_{p \in X}
V_p$ and an element $g: (v_1, V_1) \rightarrow (v_2,V_2)$ (an
arrow) in $\RR$ with $v_i \in V_i, i = 1,2$, will be a equivalence
class of triples $g= [\lambda_1,w,\lambda_2]$ where $w \in W$ for
another uniformizing system $(W,H,\rho)$, and the $\lambda_i$'s
are injections $(\lambda_i, \phi_i) : (W,H,\rho) \rightarrow
(V_i,G_i, \pi_i)$ with $\lambda_i(w) =v_i, i = 1,2$.

For another injection $(\gamma, \psi) : (W',H',\rho') \rightarrow
(W,H,\rho)$ and $w' \in W'$ with $\gamma(w')=w$ then
$[\lambda_1,w,\lambda_2] = [\lambda_1 \circ \gamma,w',\lambda_2
\circ \gamma]$

Now the maps $s,t,e,i,m$ are naturally described:
$$s([\lambda_1,w,\lambda_2]) = (\lambda_1(w),V_1), \ \ \ \ \  t([\lambda_1,w,\lambda_2]) = (\lambda_2 (w),V_2) \ \ \ \
e(x,V) = [id_V,x,id_V]$$
$$i([\lambda_1,w,\lambda_2]) =[\lambda_2,w,\lambda_1] \ \ \ \ m([[\lambda_1,w,\lambda_2],[\mu_1,z,\mu_2]) =
[\lambda_1 \circ \nu_1 , y, \mu_2 \circ \nu_2]$$
where $h=[\nu_1,y,\nu_2]$ is an arrow joining $w$ and $z$ (i.e. $\nu_1(y) =w \ \& \ \nu_2(y) =z$)

It can  be given a topology to $\RR$ so  that $s,t$ will be \'etale maps, making it
into a proper, \'etale, differentiable groupoid,  and it is not hard to check that all
the properties of groupoid are satisfied.

This is a good place to note that an orbifold $X$ given by a
groupoid $\RR\rightarrow\UU$ will be a smooth manifold if and only
if the map $(s,t)\colon\RR\rightarrow\UU\times\UU$ is one-to-one.

We  want the morphism between orbifolds to be morphisms of
groupoids, and this is precisely the case for the good maps given
in \cite{CR2}.

\vskip 0.1in \noindent {\bf Proposition 2.6 \cite{LU1}: }{\it
 A morphism of groupoids induces a good map between the underlying
orbifolds, and conversely, every good map arises in this way.}

\vskip 0.1in
\subsection{Orbispace}
    This section was taken from \cite{C}. Let $U$ be a connected, locally connected topological space.
A $G$-structure on $U$ is a triple $(\tilde{U}, G_U, \pi_U)$ where
$(\tilde{U}, G_U)$ is a connected, locally connected $G$-space and
$\pi_U: \tilde{U}\rightarrow U$ is a continuous map inducing a
homeomorphism between the orbit space $\tilde{U}/G_U$ and $U$. An
isomorphism between two $G$-structures on $U$, $(\tilde{U}_i,
G_{i, U}, \pi_{i,U})$ for $i=1,2$ is a pair $(\phi, \lambda)$
where $\lambda: G_{1,U}\rightarrow G_{2,U}$ is an isomorphism and
$\phi: \tilde{U}_1 \rightarrow \tilde{U}_2$ is a
$\lambda$-equivariant homeomorphism such that $\pi_{2,U}\circ
\phi=\pi_{1, U}$. Note that each $g\in G_U$ induces an
automorphism $(\phi_g,\lambda_g)$ on $(\widehat{U},G_U,\pi_U)$,
defined by setting $\phi_g(x)= g\cdot x,\forall x\in\widehat{U}$
and $\lambda_g(h)=ghg^{-1},\forall h\in G_U$. However, it might
not be true that every automorphism arises in this way. We shall
only take into consideration the automorphisms
$(\phi_g,\lambda_g)$, $g\in G_U$. More precisely, we {\it define}
the automorphism group of the $G$-structure
$(\widehat{U},G_U,\pi_U)$ to be $G_U$ via the induced isomorphisms
$(\phi_g,\lambda_g)$ on it. We would like to point out that since
the action of $G_U$ on $\widehat{U}$ is not required to be
effective, two {\it different} automorphisms of the $G$-structure
could have the {\it same} induced map.

\vspace{1.5mm}

Given a $G$-structure $(\widehat{U},G_U,\pi_U)$ on $U$, we
consider the inverse image $\pi_U^{-1}(W)$ in $\widehat{U}$, where
$W$ is a connected open subset of $U$. Denote by $\widehat{W}$ one
of the connected components of $\pi_U^{-1}(W)$, by $G_W$ the
subgroup of $G_U$ consisting of elements $g\in G_U$ such that
$g\cdot\widehat{W}=\widehat{W}$, and let
$\pi_W=(\pi_U)|_{\widehat{W}}$.

\vspace{1.5mm}

\noindent{\bf Lemma 2.7:} {\em The triple
$(\widehat{W},G_W,\pi_W)$ defines a $G$-structure on $W$.
Moreover, the action of $G_U$ on $\widehat{U}$ induces a
transitive action on the set of all such $G$-structures on $W$ for
which the following holds: Let $(\widehat{W_i}, G_{W,i},
\pi_{W,i})$, $i=1,2$, be two such $G$-structures, and
$\widehat{W_2}=g\cdot\widehat{W_1}$ for some $g\in G_U$, then
$G_{W,2}=gG_{W,1}g^{-1}$ in $G_U$. The stabilizer of the action at
$G$-structure $(\widehat{W},G_W,\pi_W)$ is precisely the subgroup
$G_W$ in $G_U$. }

\vspace{2mm}

We will say that $(\widehat{W},G_W,\pi_W)$ is {\it induced} from
the $G$-structure $(\widehat{U},G_U,\pi_U)$. Note that the
subgroup $G_W$ is both closed and open in $G_U$. In fact, the
space of cosets $G_U/G_W=\{gG_W|g\in G_U\}$ inherits a discrete
topology from $G_U$.

\vspace{1.5mm}

{\it Transition maps.}\hspace{2mm} Let $U_1$, $U_2$ be connected
open subsets of a locally connected topological space with
$U_1\cap U_2\neq \emptyset$. Suppose $U_1$ and $U_2$ are given
with $G$-structures $(\widehat{U_1},G_{U_1},\pi_{U_1})$ and

$(\widehat{U_2},G_{U_2},\pi_{U_2})$ respectively. Let $\{W_i|i\in
I\}$ be the set of connected components of $U_1\cap U_2$. We
define $Iso_{W_i}(U_1,U_2)$ to be the set of isomorphisms from an
element of $\G_{U_1}(W_i)$ to an element of $\G_{U_2}(W_i)$, and
set
$$Iso(U_1,U_2)=\bigsqcup_{i\in I} Iso_{W_i}(U_1,U_2). \leqno (2.1)
$$
There is a canonical right action of $G_{U_1}\times G_{U_2}$ on
$Iso(U_1,U_2)$ defined by
$$
((\phi,\lambda),(g_1,g_2))\mapsto (\phi,\lambda)\cdot (g_1,g_2):=
(g_2^{-1}\circ\phi\circ g_1,Ad(g_2^{-1})\circ\lambda\circ
Ad(g_1)).
$$
Note that $Domain(\lambda)$ is contained in the stabilizer at
$(\phi,\lambda)$, through the identification to a subgroup of
$G_{U_1}\times G_{U_2}$ by $g\mapsto (g,\lambda(g))$.

\vspace{1.5mm}

\noindent{\bf Definition 2.8:} {\em {\em A set of transition maps
from $U_1$ to $U_2$} is a topological space $Tran(U_1,U_2)$
together with a map
$$
\chi_{U_1U_2}:Tran(U_1,U_2)\rightarrow Iso(U_1,U_2) \leqno (2.2)
$$
satisfying the following conditions:
\begin{itemize}
\item [{(a)}]
For any $i\in I$,
$Tran_{W_i}(U_1,U_2):=\chi_{U_1U_2}^{-1}(Iso_{W_i}(U_1,U_2))$ is a
non-empty open subset.
\item [{(b)}] There is a continuous right action of $G_{U_1}\times G_{U_2}$
on $Tran(U_1,U_2)$, written $(\xi,(g_1,g_2))\mapsto \xi\cdot
(g_1,g_2)$, which is transitive when restricted to each
$Tran_{W_i}(U_1,U_2)$, such that the map $\chi_{U_1U_2}$ is
$G_{U_1}\times G_{U_2}$-equivariant.
\item [{(c)}] The stabilizer at $\xi\in Tran(U_1,U_2)$, $Stab(\xi)$,
is the image of the injective homomorphism
$Domain(\lambda)\rightarrow G_{U_1}\times G_{U_2}$ defined by
$g\mapsto (g,\lambda(g))$ where $\lambda$ is given in
$\chi_{U_1U_2}(\xi)=(\phi,\lambda)$.
\item [{(d)}] For any $i\in I$ and $\xi\in Tran_{W_i}(U_1,U_2)$, the map
$G_{U_1}\times G_{U_2}\rightarrow Tran_{W_i}(U_1,U_2)$ defined by
$(g_1,g_2)\mapsto \xi\cdot (g_1,g_2)$ induces a homeomorphism
between the space of left cosets $(G_{U_1}\times
G_{U_2})/Stab(\xi)$ and $Tran_{W_i}(U_1,U_2)$.
\end{itemize}
Each element of $Tran(U_1,U_2)$ is called a {\em transition map}.
}

\hfill $\Box$

As an example, let us consider the case when $U_1=U_2=U$. Define a
map $\chi_U: G_U\rightarrow Iso(U,U)$ by
$$
\chi_U(g)=(g,Ad(g)), \leqno (2.3)
$$
and a right action of $G_U\times G_U$ on $G_U$ by $g\cdot
(g_1,g_2) =g_2^{-1}gg_1$. Then one can easily verify that
$(G_U,\chi_U)$ satisfies Definition 2.8, hence defines a set of
transition maps from $U$ to $U$.

\vspace{1.5mm}

The {\it composition} of transition maps is a map
$Tran(U_1,U_2)\times Tran(U_2,U_3)\rightarrow Tran(U_1,U_3)$,
written $(\xi_{12},\xi_{23})\mapsto \xi_{23}\circ\xi_{12}$, such
that $\chi_{U_2U_3}(\xi_{23})\circ \chi_{U_1U_2}(\xi_{12})$ equals

$\chi_{U_1U_3}(\xi_{23}\circ\xi_{12})$ when restricted to the
intersection of their domains. However, this could be problematic
in general. First of all, in order for
$\chi_{U_2U_3}(\xi_{23})\circ \chi_{U_1U_2}(\xi_{12})$ to be
defined, the following condition
$$
Range(\chi_{U_1U_2}(\xi_{12}))\cap
Domain(\chi_{U_2U_3}(\xi_{23}))\neq \emptyset \leqno (2.4)
$$
is necessary, which does not hold in general. Secondly, the
inverse image of the left-hand side of $(2.4)$ under
$\chi_{U_1U_2}(\xi_{12})$, which is a subset of $\G_{U_1}(U_1\cap
U_2\cap U_3)$, could be contained in more than one element of
$\G_{U_1}(U_1\cap U_3)$. In this case, the composition map will
pick up one of these elements of $\G_{U_1}(U_1\cap U_3)$ for the
domain of $\chi_{U_1U_3}(\xi_{23}\circ\xi_{12})$, which seems
quite unnatural.

Because of these reasons, we shall impose composition maps
$Tran(U_1,U_2)\times Tran(U_2,U_3)\rightarrow Tran(U_1,U_3)$ only
when the following condition for $U_1,U_2,U_3$ is met
$$
\begin{array}{ccccc}
(1) & & & U_1\subset U_2\subset U_3 & \mbox {or}\\
(2) & & & U_3\subset U_2\subset U_1 & \mbox {or}\\
(3) & & & U_2\subset U_1, U_2\subset U_3. & \\
\end{array}
\leqno (2.6)
$$
It turns out that $(2.6)$ is sufficiently rich so that the
composition of transition maps in general can be derived from it,
but on the other hand, $(2.6)$ is also minimal that none of the
cases in $(2.6)$ may be removed.

\vspace{1.5mm}

 With these
preparations, we now introduce

\vspace{1.5mm}

\noindent{\bf Definition 2.9:} {\em Let $X$ be a locally connected
topological space. An {\em orbispace structure} on $X$ is a set
$\U$ of connected open subsets of $X$ satisfying the following
conditions:
\begin{enumerate}
\item For any $U_1,U_2\in\U$ with $U_1\cap U_2\neq \emptyset$, each
connected component of $U_1\cap U_2$ is also in $\U$.
\item Each element $U$ of $\U$ is assigned with a $G$-structure $(\widehat{U},
G_U,\pi_U)$ satisfying the following conditions:
\begin{itemize}
\item [{a)}] Any ordered pair of $U_1,U_2\in\U$ with
$U_1\cap U_2 \neq\emptyset$ is assigned with a set of transition
maps $(Tran(U_1,U_2),\chi_{U_1 U_2})$ as defined in Definition
2.8.
\item [{b)}] For any $U_1,U_2,U_3\in \U$ satisfying $(2.6)$, there
is a composition map $Tran(U_1,U_2)\times Tran(U_2,U_3)\rightarrow
Tran(U_1,U_3)$, written $(\xi_{12},\xi_{23})\mapsto
\xi_{23}\circ\xi_{12}$, such that
$\chi_{U_1U_3}(\xi_{23}\circ\xi_{12})$ equals
$\chi_{U_2U_3}(\xi_{23})\circ \chi_{U_1U_2}(\xi_{12})$  when
restricted to the domain of the latter. Moreover, the following
condition is satisfied: for any $g_1\in G_{U_1}$,
$g_2,g_2^\prime\in G_{U_2}$, and $g_3\in G_{U_3}$,
\begin{eqnarray*}
(\xi_{23}\cdot (g_2^\prime,g_3))\circ (\xi_{12}\cdot (g_1,g_2)) &
= & ((\xi_{23}\cdot (g_2^\prime
g_2^{-1},1_{G_{U_3}}))\circ\xi_{12})
\cdot (g_1,g_3) \\
& = & (\xi_{23}\circ (\xi_{12}\cdot
(1_{G_{U_1}},g_2(g_2^\prime)^{-1}))) \cdot (g_1,g_3).
\end{eqnarray*}
Combined with Definition 2.8 (d), this condition implies that the
composition map is continuous.
\item [{c)}] The composition map in b) is associative in the following
sense: for any $U_1,U_2,U_3,U_4\in\U$ such that $(U_1,U_2,U_3)$,
$(U_1,U_3,U_4)$, $(U_1,U_2,U_4)$ and $(U_2,U_3,U_4)$ satisfy
(2.1.6),
$$
\xi_{34}\circ
(\xi_{23}\circ\xi_{12})=(\xi_{34}\circ\xi_{23})\circ\xi_{12}.
$$
\item [{d)}] For any $U\in\U$, $(Tran(U,U),\chi_{UU})=(G_U,\chi_U)$
(cf. (2.4)). Moreover, the composition map in b) coincides with
the multiplication in $G_U$, i.e., for any $g_1,g_2\in G_U$,
$$
g_1\circ g_2=g_1g_2.
$$
We shall call an element of $Tran(U,U)=G_U$ an {\em automorphism}
of $(\widehat{U},G_U,\pi_U)$.
\end{itemize}
\item The set $\U$ is a cover of $X$.
\end{enumerate}

The topological space $X$ equipped with the orbispace structure
$\U$ is called an {\em orbispace}, which will be denoted by
$(X,\U)$ in general. }

\hfill $\Box$

Each element $U$ of $\U$ is called a {\it basic open set} of $X$.
For any point $p\in X$, which is contained in a basic open set
$U$, pick a $x\in \widehat{U}$ in the inverse image
$\pi_U^{-1}(p)$. We define the {\it local group} of $p$ to be the
stabilizer $G_x$ of $x$ in $G_U$ (i.e.  $G_x=\{g\in G_U|g\cdot
x=x\}$), and denote it by $G_p$. Clearly different choices of $x$
result in the same conjugacy class in $G_U$, and different choices
of $U$ give rise to isomorphic groups because of the existence of
transition maps. An orbispace structure $\U$ is called {\it
trivial} if $G_U$ is trivial for each $U\in\U$. (Every locally
connected topological space is canonically an orbispace with a
trivial orbispace structure.) An orbispace $(X,\U)$ is called {\it
orbifold} if $G_U$ is finite for each $U\in\U$. As a notational
convention, we very often only write $X$ for an orbispace
$(X,\U)$, and write $X_{top}$ for the underlying topological space
$X$ for simplicity.

    \vskip 0.1in \noindent{\bf Example 2.10:}{\it For any locally
connected $G$-space $(Y,G)$, the orbit space $Y/G$ canonically
inherits an orbispace structure, of which $\U$ is taken to be the
set of all connected open subsets of the orbit space $Y/G$. The
$G$-structure assigned to each element of $\U$ is a fixed choice
of the $G$-structures induced from the God-given $G$-structure
$(Y,G)$ on the orbit space $Y/G$. Each set of transition maps is
obtained by restricting $G$ to the corresponding induced
$G$-structures. The verification of Definition 2.8 for this case
is straightforward. An orbispace will be called {\it global} if it
arises as the orbit space of a $G$-space equipped with the
canonical orbispace structure as discussed in this example.}

\vspace{1.5mm}

\noindent{\bf Remark 2.11: }{\it When a basic open set $U_\alpha$
is included in another one $U_\beta$, a transition map in
$Tran(U_\alpha,U_\beta)$ is just an isomorphism from the
$G$-structure of $U_\alpha$ onto one of the $G$-structures of
$U_\alpha$ induced from the $G$-structure of $U_\beta$. We will
call a transition map in $Tran(U_\alpha,U_\beta)$ an {\it
injection} when $U_\alpha\subset U_\beta$ holds.}

\vskip 0.1in

\noindent{\bf Remark 2.12: }{\it In order to minimize the
dependence of orbispace structure on the choice of the set $\U$ of
basic open sets so that it becomes more intrinsic, it is
appropiate to introduce an equivalence relation as follows: Two
orbispace structures are said to be {\it directly equivalent} if
one is contained in the other. The equivalence relation to be
introduced is just a finite chain of direct equivalence.}

\vskip 0.1in \noindent {\bf Remark 2.13: }{\it The definition of
good map/morphism  can be modified from that of orbifold
\cite{CR2} in an obvious fashion. } \vskip 0.1in

In a certain sense, the definition of orbispace structure in this
paper is of a less intrinsic style than the definition of orbifold
structure originally used in \cite{S}, where the language of germs
was used. The reason for which we choose such a formulation lies
in the following facts: First of all, an automorphism of a
$G$-structure $(\widehat{U},G_U,\pi_U)$, i.e., a pair
$(\phi,\lambda)$ where $\lambda\in Aut(G_U)$ and $\phi$ is a
$\lambda$-equivariant homeomorphism of $\widehat{U}$ inducing
identity map on $U$, might not be induced by an action of $G_U$ on
$\widehat{U}$ (at least we have a burden of proving so). But we
want to only consider those arising from the action of $G_U$.
Secondly, suppose $W$ is a connected open subset of $U$, both of
which are given with $G$-structures $(\widehat{W},G_W,\pi_W)$ and
$(\widehat{U},G_U,\pi_U)$ respectively. Then an open embedding
$\phi:\widehat{W}\rightarrow \widehat{U}$, which is
$\lambda$-equivariant for some monomorphism
$\lambda:G_W\rightarrow G_U$ and induces the inclusion
$W\hookrightarrow U$, might not be an isomorphism onto one of the
$G$-structures of $W$ induced from $(\widehat{U},G_U,\pi_U)$.
However, this kind of pathology can be ruled out if each
$\widehat{U}$ is Hausdorff, $G_U$ acts on $\widehat{U}$
effectively. \vspace{1.5mm}

\section{Cohomology ring of crepant resolutions}

    Suppose that $X$ is an orbifold. In general, $K_X$ is an orbifold vector bundle or a $Q$-divisor only.
    When the $X$ is so called  Gorenstein, $K_X$ is a bundle or a divisor. For the Gorenstein orbifold, a
    resolution $\pi:Y \rightarrow X$ is called a crepant resolution if $\pi^*K_X=K_Y$. Here, "crepant"
    can be viewed as a minimality condition with respect to the canonical bundle. A crepant resolution always
    exists when dimension is two or three. A nice way to construct it is to use Hilbert scheme of points.
    However, the crepant resolution in dimension three is not unique. Different crepant resolutions are
    connected by flops. When the dimension is bigger than four, the crepant resolution does not always exist.
    It is an extremely interesting problem in algebraic geometry to find out when it does exist.

When the orbifold string theory was first constructed over the global
quotient \cite{DHVW}, one of its first invariant is
orbifold Euler characteristic. It was conjectured that the
orbifold Euler characteristic is the same as the Euler
characteristic of its crepant resolution. This fits well with
McKay's correspondence in algebro-geometry. It had been the main
attraction before the current development. By the work of Batyrev
and others \cite{B1} \cite{DL}, this conjecture has been extended
and solved for the orbifold Hodge number of Gorenstein global
quotients. Very recently, it was solved in the complete
generalities by Lupercio-Poddar \cite{LP} and Yasuda \cite{Y}.
Among all the examples, Hilbert scheme of points of algebraic
surfaces are particularly interesting.
  Suppose that $M$ is an algebraic surface. We use $M^{[n]}$ to denote the Hilbert scheme of points of length
    $n$ of $M$. In his thesis \cite{G}, G\"{o}ttsche computed the generating function of the Euler number $\sum_{n=1}^{\infty} \chi(M^{[n]})q^n$
and showed that it has a surprising modularity. In 1994, in order
to explain its modularity, Vafa-Witten \cite{VW} computed
$\H=\oplus_n H^*(M^n/S_n,\C)$. Motivated by the 
    orbifold conformal field theory, they directly wrote $\H$ as a "Fock space" or a representation
    of the Heisenberg algebra. Then, the generating function of the Euler characteristic is interpreted as the correlation function of
    an elliptic curve. Therefore, it should be invariant under the modular transformations of the elliptic curve. This shows that
the space of cohomology itself has more structures. The orbifold string
    theory conjecture predicates that $\oplus_n H^*(M^{[n]}, \C)$ should also admit a representation of the Heisenberg algebra.
    This conjecture was verified by a beautiful work of Nakajima \cite{N} and others.

    The ring structure of $M^{[n]}$ is quite subtle and more
    interesting. Partial results have been obtained by
    Fantechi-G\"{o}ttsche \cite{FG1}, Ellingsrud-Stromme
    \cite{ES1} \cite{ES2}, Beauvill \cite{Bea}, Mark \cite{Mar}. Based on  an important observation by Frenkel
    and Wang \cite{FW}, Lehn-Sorger \cite{LS1} determined the
    cohomology ring of $(\C^2)^{[n]}$. At the same time, the
    author was computing orbifold cohomology $(\C^2)^n/S_n$ and
    the result  from both calculations match perfectly.  Based on the physical motivation and the strong evidence from $(\C^2)^n/S_n$.
    the author proposed \cite{R2} a conjecture in the case
    of a hyperkahler resolution.
    \vskip 0.1in
    \noindent
    {\bf Cohomological Hyperkahler Resolution Conjecture: }{\it Suppose that $\pi: Y\rightarrow X$ is a hyperkahler resolution. Then, the ordinary cohomology ring of $Y$
    is isomorphic to the orbifold cohomology ring of $X$.}
    \vskip 0.1in
    In the case of the Hilbert scheme points of surfaces, the Cohomological Hyperkahler Resolution Conjecture (CHRC) implies that $K3^{[n]}, (T^4)^{[n]}$ have an isomorphic cohomology ring
    as the orbifold cohomology rings of $K3^n/S_n, (T^4)^{n}/S_n$. The later was proved recently by the beautiful works
    of Lehn-Sorge \cite{LS2}, Fantechi-G\"{o}ttsche \cite{FG2} and Uribe \cite{U}.
    It should be mentioned that  Fantechi-G\"{o}ttsche-Uribe's work computed the orbifold ring structure of $X^n/S_n$ for an arbitrary
    complex manifold $X$ which may or may not be $K3, T^4$. There is a curious phenomenon that over rational number Lehn-Sorge, Fantanch-G\"{o}ttsche and Uribe showed that one must modify the ring structure of the orbifold
    cohomology by a sign in order to match to the cohomology of Hilbert scheme. However, Qin-Wang observed that
	such a sign modification is unnecessary over complex number \cite{QW}. All the conjectures stated in this article (in fact
	any conjecture motivated by physics) are the statements over complex number. 

 In a different direction, the Hilbert scheme can
    be regarded as a moduli space. For a moduli space, there is a conjectural {\em Mumford principal} that
    the K\"{u}nneth components of the Chern characters of the universal sheaf form a set of ring generators.
    Such a Mumford principal was well-known in the Donaldson theory. For the Hilbert scheme of points, the Mumford principal
    was proved by Li-Qin-Wang \cite{LQW1} which is instrumental  in Lehn-Sorger's calculation of $M^{[n]}$ for
    $M=K3, T^4$ (see \cite{LQW2} for another set of generators). However, it is quite difficult to determine the relations
    for Li-Qin-Wang generators. The ring structure of $X^{[n]}$ for a general algebraic
    surface $X$ is     still unknown.

    It is easy to check that CHRC is false if we drop the hyperkahler condition. Motivated by physics and the work of the author with An-Min Li
    on the quantum cohomology and flop, the author proposed a conjecture for the arbitrary crepant resolution.

    As  mentioned previously, the crepant resolutions are not unique. The different crepant resolutions are connected by
    "K-equivalence" \cite{W}.  Two smooth (or Gorenstein orbifolds) complex
    manifolds $X, Y$ are $K$-equivalent iff there is a common resolution $\phi, \psi: Z\rightarrow X, Y$ such that $\phi^* K_X=
    \psi^* K_Y$. Batyrev-Wang \cite{B}, \cite{W} showed that two $K$-equivalent projective manifolds have the same betti number. It is natural to
    ask if they have the same ring structures. This question is obviously related to CHRC. Suppose that CHRC holds for non-hyperhahler
    resolutions. It implies that different resolutions (K-equivalent) have the same ring structures. Unfortunately, they usually have
    different ring structures,  and hence CHRC fails in general. It is
easy to check this
    in case of  three dimensional flops. A key idea to remedy the  situation is to include the quantum corrections. The author proposed \cite{R1}
    \vskip 0.1in
    {\bf Quantum Minimal Model Conjecture: }{\it Two $K$-equivalent projective manifolds have the same quantum cohomologies.}
    \vskip 0.1in
    Li and the author proved Quantum Minimal Model conjecture in complex dimension three. In higher dimensions, it seems to be
    a difficult problem. In many ways, Quantum Minimal Model Conjecture unveils the deep relation between the quantum cohomology and
    the birational geometry \cite{R1}. However, it is a formidable task to master the quantum cohomology machinery for any non-experts.
     In \cite{R4}, the author proposed another conjecture focusing on the cohomology
    instead of the quantum cohomology. As mentioned before, the cohomology ring structures are not isomorphic for $K$-equivalent
    manifolds. Therefore, some quantum information must be included. The new conjecture requires a minimal set of quantum information
    involving the GW-invariants of the exceptional
    rational curve.

\subsection{Orbifold cohomology of symmetry product}

    The orbifold cohomology of symmetry product is a beautiful subject. In this section, we first present Vafa-Witten's calculation of the orbifold cohomology 
 group $\H=\oplus_n H^*_{orb}(X^n/S_n, \R)$ and directly write it as
a "Fock space" or the representation of super Heisenberg algebra.
Then, we sketch Fantechi-G\"{o}ttsche \cite{FG2} and Uribe
\cite{U}'s calculation of the ring structure.

    Let $X$ be an almost complex manifold and $S_n$ be the symmetry group on $n$-letters. Then, $X^n/S_n$ is an almost complex
orbifold and its orbifold cohomology $H^*_{orb}(X^n/S_n, \R )$ is well-defined.
    \vskip 0.1in
    \noindent
    {\bf Theorem 3.1: }{\it $\H=\oplus_n H^*_{orb}(X^n/S_n,
\R)$ is an irreducible representation of the super Heisenberg algebra.
}
    \vskip 0.1in
    \noindent
    {\bf Proof: }
     We first compute the cohomology of nontwisted sector
$H^*(X^n/S_n, \R)$. It is easy to see that $H^*(X^n/S_n,
\R)=H^*(X^n, \R)^{S_n}$. Pick a basis $w^a$ of the cohomology of
$X$. Then, $w^{a_1}\otimes w^{a_2}\dots\otimes w^{a_n}\in H^*(X^n,
\R)$. We obtain a class of $H^*(X^n, \R)^{S_n}$ by symmetrizing
$w^{a_1}\otimes w^{a_2}\dots\otimes w^{a_n}$. We denote this class
using symbol $\alpha^{a_1}_{-1}\cdots \alpha^{a_n}_{-1}|0>$. This
is a physical notation. $0>$ is called {\em vacum} representing
the cohomology of $H^*(X^0/S_0, \R)$ and $\alpha^{a}_{-1}$ is
called {\em a creation operator}. We think formally that
$\alpha^a_{-1}$ acts on vacum $|0>$ to create  an element of
$H^*_{orb}(X^n/S_n, \R)$ called {\em 1-particle state}. The
creation operators satisfy commuatation relation
$\alpha^a_{-1}\alpha^b_{-1}=(-1)^{deg(a)deg(b)}\alpha^b_{-1}
\alpha^a_{-1}$. We will see this commutation relation again when
we compute the cohomology of twisted sector.

    The twisted sector is given by the connected components of $X^n_g/C(g)$ for the conjugacy class $(g)$. It is well-known that
a conjugacy class is uniquely determined by the cycle
decomposition $g=1^{n_1}2^{n_2}\cdots k^{n_k}$ where $i^{n_i}$
means the $n_i$-many cycles of length $i$. It is clear that
$n=\sum_i i n_i$. The fixed point loci $X^n_g=X^{n_1}\times \cdots
\times X^{n_k}$. The centralizer $C(g)=S_{n_1}\bar{\times}
Z_1\cdots \times S_{n_k}\bar{\times} Z_{n_k}$, where
$\bar{\times}$ means the semi-direct product. Hence, as a
topological space, the twisted sector $X^n_{(g)}=X^n_g/
C(g)=X^{n_1}/S_{n_1}\cdots\times X^{n_k}/S_{n_k}$ although it has
different orbifold structure given by the extra group $Z_1, \cdots
Z_{n_k}$. The cohomology of the factor $X^{n_i}/S_{n_i}$ is called
{\em $k$-particle state}. It can be written down by the same
method of the nontwisted sector. We use a lower indice $\alpha^a_{-i}$
to denote it. A cohomology class of a twisted sector can be written
as
$$\alpha^{a^1_1}_{-1}\cdots \alpha^{a^1_{n_1}}_{-1}\alpha^{a^2_1}_{-2}\cdots \alpha^{a^k_{n_k}}_{-k}|0>.\leqno(3.1)$$
When we consider the direct sum $\H=\oplus_n H^*_{orb}(X^n/S_n,
\R)$, the constraint $n=\sum_i n_i$ becomes a trivial condition.
Hence, we conclude that the elements of the form (3.1) forms a
basis of $\H$. If you are a physicist, you stop here and claim
$\H$ is a "Fock space" or an irreducible representation of super
Heisenberg algebra. Since a common knowledge in physics may not be
a common knoweldge in mathematics, let me reconstruct the action
of the Heisenberg algebra which is motivated by harmonic oscillator.

    Before we reconstruct the action of super Heisenberg algebra, we compute the degree shifting number. Then, we will see that the theory
is different when $\dim X$ is even or odd. For each $j$-cycle, its
action on $(C^N)^j$ has eigenvalues $N$-copy of $e^{\frac{2\pi i
p}{j}}$, for $p=0, \cdots, j$. Therefore, its contribution to
degree shifting number is $\frac{j-1}{2}N$. Let $deg (g)$ be the
minimal number of transpositions to express $g$ as the composition of
transpositions. It is clear that $deg(g)=\sum_i n_i(i-1)$ and
degree shifting number $\iota_{ (g)}=\frac{N}{2} deg(g)$. Notes
that when $N$ is even, $\iota_{(g)}$ is an integer. Otherwise, it
is a fraction. In particular, when $N=2$, $\iota_{(g)}=deg(g)$.

    To construct the action of the full super Heisenberg algebra, we
    need to add so called {\em annihilation operators}. Let $H$ be
    a superspace, i.e., $H=H_{even}\oplus H_{odd}$. We further assume that $H$ is
    equipped with an inner product such as Poincare paring. For any $a\in H_{even}$,
    we define $deg(a)=0$ and $deg(b)=1$ if $b\in H_{odd}$. A super
    Heisenberg algebra is the set of operators $\alpha^a_{l}$ for
    $a\in H, l\in \Z-\{0\}$. We define super commutator
    $$\{\alpha^a_l, \alpha^b_m\}=\alpha^a_l\alpha^b_m -
    (-1)^{deg(a)deg(b)}\alpha^b_m \alpha^a_l.\leqno(3.2)$$
    Then, the operators $\alpha^a_l$ satisfies relation
    $$\{\alpha^a_l, \alpha^b_m\}=-l\delta_{l+m, 0} <a,b>Id.\leqno(3.3)$$
    The operators $\alpha^a_l$ for $l>0$ are called annihilation
    operators and $\alpha^a_l$ for $l<0$ are called creation
    operators. The action of the annihilation operators on $\H$ is
    determined by the commutation relation of (3.3) and
    $$\alpha^a_l|0>=0$$
    for $l>0$. To recover the Euler characteristic, we define
    operator $L_0$ by the condition
    $$L_0|0>=1, \{L_0, \alpha^a_l\}=-l \alpha^a_l.\leqno(3.4)$$
    Its character
    $$tr q^{L_0}=\sum q^n (\dim (n-\mbox{even eigenspace})-\dim
    (n-\mbox{odd eigenspace})).\leqno(3.5)$$
    From the general theory of Heisenberg representation, $tr
    q^{L_0}$ is a modular form.

    It is easy to check that
    $$L_0 (\alpha^{a_1}_{-l_1}\cdots
    \alpha^{a_k}_{-l_k}|0>)=\sum_i l_i \alpha^{a_1}_{-l_1}\cdots
    \alpha^{a_k}_{-l_k}|0>.$$
    Let $n=\sum_i$. Then,  $\alpha^{a_1}_{-l_1}\cdots
    \alpha^{a_k}_{-l_k}|0>\in H^*_{orb}(X^n/S_n, \R)$.
    Furthermore, when the dimension of $X$ is even, the degree
    shifting does not change the even/odd property of orbifold
    cohomology classes for its original degree. Hence,
    $$q^{L_0}=\sum_n q^n \chi (H^*_{orb}(X^n/S_n, \R).\leqno(3.6)$$
    An routine calculation shows that it equals to
    $$\prod_n \frac{1}{(1-q^n)^{\chi(X)}}.\leqno(3.7)$$

    Next, we calculate its ring structure.  Although their goal is to compute the cohomology of
    the Hilbert scheme of points, the algebraic structure  set up by Lehn-Sorger \cite{LS2} is very
    convenient for the computation of the orbifold cohomology of a symmetry product. The
    actual computation of the orbifold product was carried out independently by
    Fantechi-G\"{o}ttsche \cite{FG2} and Uribe \cite{U}.

    Before we carry out the
    computation, we take the advantage of the global quotient to
    reformulate the orbifold cohomology slightly.

    Let $Y=X/G$ be an orbifold with $X$ a compact
complex manifold. As before $X^g$ will denote the fixed point set
of the action of $g$ on $X$.

The cohomology classes will be labeled by elements in $G$ and let
the total ring $A(X,G)$ be
$$A(X,G) := \bigoplus_{g \in G} H^*(X^g;\C) \times \{g\}$$
Its group structure is the natural one and the ring structure that
will be defined later will give us the orbifold cup product. The
grading is the one  in the orbifold cohomology, i.e.
$$A^d(X,G) = \bigoplus_{g \in G} H^{d-2\iota_{(g)}}(X^g;\C) \times \{g\}$$

For $h \in G$ there is a natural map $h: X^g \rightarrow
X^{hgh^{-1}}$ which can be extended to an action in $A(X,G)$
inducing an isomorphism

\begin{eqnarray*} h :H^*(X^g;\C) \times \{g\} &
\rightarrow & H^*(X^{hgh^{-1}};\C) \times \{hgh^{-1}\}\\
(\alpha,g) & \mapsto & ((h^{-1})^*\alpha, hgh^{-1})
\end{eqnarray*}

The invariant part under the action of $G$ is isomorphic as a
group to the orbifold cohomology,

$$A(X,G)^G \cong \bigoplus_{(g)} H^*(X^g;\C)^{C(g)} \cong
H^*_{orb}(X/G;\C)$$

The construction of orbifold cup product use quotient $X^g/C(g)$.
We can obviously lift it to the fixed point set.

For $X^{\langle h_1,h_2 \rangle}$, the fixed point set of $\langle
h_1,h_2 \rangle$,
 let

 $$f^{h_i,\langle h_1,h_2 \rangle} :
H^*(X^{h_i};\C) \rightarrow H^*(X^{\langle h_1, h_2 \rangle};\C)$$
$$f_{\langle h_1,h_2 \rangle, h_i} :
H^*(X^{\langle h_1, h_2 \rangle};\C) \rightarrow H^*(X^{h_i};\C)$$
be the pull-back and the push-forward respectively of the diagonal
inclusion map $X^{\langle h_1, h_2 \rangle} \hookrightarrow
X^{h_i}$ where $i=1,2,3$ and $h_3=(h_1h_2)$.

We need to make use of the obstruction bundle over $X^{\langle
h_1, h_2 \rangle}$;
 as $Y_{(h_1,h_2)} = X^{\langle h_1, h_2 \rangle}/C(h_1,h_2)$ and taking the
 projection map $\pi :X^{\langle h_1, h_2 \rangle} \rightarrow X^{\langle h_1, h_2 \rangle}/C(h_1,h_2)$ we
 will consider the Euler class of $\pi^*(E_{({\bf h})})$.

Let the product $A(X,G) \otimes A(X,G)
\stackrel{\cdot}{\rightarrow} A(X,G)$ be defined by
$$(\alpha,h_1) \cdot (\beta,h_2) := f_{\langle h_1,h_2 \rangle, h_1h_2} \left(
f^{h_1,\langle h_1,h_2 \rangle} (\alpha) \wedge f^{h_2,\langle
h_1,h_2 \rangle}(\beta)  \wedge \pi^*c(E_{({\bf h})}) \right)$$
whose three point function is
\begin{eqnarray*}
\lefteqn{<(\alpha,h_1),(\beta,h_2),(\gamma,(h_1h_2)^{-1})>} & &\\
& := & \int_{X^{\langle h_1, h_2 \rangle}}
 f^{h_1,\langle h_1,h_2 \rangle} (\alpha) \wedge f^{h_2,\langle h_1,h_2
\rangle}(\beta)  \wedge f^{(h_1h_2)^{-1}, \langle h_1,h_2
\rangle}(\gamma) \wedge \pi^*c(E_{({\bf h})})
\end{eqnarray*}

The product $A(X,G) \otimes A(X,G) \stackrel{\cdot}{\rightarrow}
A(X,G)$ previously defined is $G$ equivariant.

This product induces a ring structure on the invariant group
$A(X,G)^G$ which will match with the orbifold cup product. Thus
$A(X,G)^G$ will inherit the properties of the orbifold cup
product.

$X$ will be an even dimensional complex manifold
 $dim_{\C} X =2N$, and the orbifold in
mind will be $X^n/S_n$ where the action of the symmetric group
$S_n$ on $X^n$ is the natural one.

Now, we introduce more notations. For $\sigma, \rho \in S_n$, let
$\Gamma \subset [n]:=\{1,2,\dots,n\}$ be a set stable under the
action of $\sigma$; we will denote by $\O(\sigma;\Gamma)$ the set
of orbits induced by the action of $\sigma$ in $\Gamma$. If
$\Gamma$ is $\sigma$-stable and $\rho$-stable,
$\O(\sigma,\rho;\Gamma)$ will be the set of orbits induced by
$\langle \sigma,\rho \rangle$. When the set $\Gamma$ is dropped
from the expression, the set $\O(\sigma, [n])$ will be  denoted
$\O(\sigma)$.

$|\sigma|$  will denote  the minimum number $m$ of transpositions
$\tau_1, \dots ,\tau_m$ such that $\sigma=\tau_1\dots\tau_m$;
hence
$$|\sigma|+|\O(\sigma)|=n$$

The set $X^n_\sigma$ will denote the fixed point set under the
action of $\sigma$ on $X^n$. Superscripts on $X$ will count the
number of copies of itself on the cartesian product, and
subscripts will be elements of the group and will determine fixed
point sets.

For $h_1,h_2 \in S_n$ the obstruction bundle $E_{({\bf h})}$ over
$Y_{(h_1,h_2)}$ is defined by
$$E_{({\bf h})}= \left( H^1(\Sigma) \otimes e^*TY \right)^G$$
where $G= \langle h_1,h_2 \rangle$, $Y=X^n/S_n$ and $\Sigma$ is an
orbifold Riemann surface provided with a $G$ action such that
$\Sigma / G =(S^2, (x_1,x_2,x_3),(k_1,k_2,k_3))$ is an orbifold
sphere with three marked points.

Because $H^1(\Sigma)$ is a trivial bundle, the pullback of
$E_{({\bf h})}$ under $\pi: X^n_{h_1,h_2} \rightarrow
Y_{(h_1.h_2)}$ is
$$E_{h_1,h_2}:=\pi^*E_{({\bf h})} = \left( H^1(\Sigma) \otimes \Delta^*TX^n \right)^G$$
where $\Delta: X^n_{h_1,h_2} \hookrightarrow X^n$ is the inclusion
(if $\rho : X^n \rightarrow Y$ is the quotient map, then $\rho
\circ \Delta = e \circ \pi$).

Without loss of generality we can assume that $|\O(h_1,h_2)|=k$
and $n_1+ \cdots + n_k =n$ a partition of  such that
$$\Gamma_i =\{n_1 + \cdots + n_{i-1}+1, \dots , n_1 + \cdots +n_i
\}$$ and $\{ \Gamma_1, \Gamma_2, \dots , \Gamma_k \} =
\O(h_1,h_2)$. We observe that  the obstruction bundle
$E_{h_1,h_2}$ is the product of $k$ bundles over $X$ (i.e.
$E_{h_1,h_2} = \prod_i E_{h_1,h_2}^i$), where each factor $E_{h_1,
h_2}^i$ corresponds to the orbit $\Gamma_i$.

For $\Delta_i : X \rightarrow X^{n_i}$ $i=1,\dots,k$  the diagonal
inclusions, the bundles $\Delta_i^*TX^{n_i}$ become $G$ bundles
via the restriction of the action of $G$ into the orbit $\Gamma_i$
and
$$\Delta^*TX^n \cong \Delta_1^*TX^{n_1} \times \cdots \times
\Delta_k^*TX^{n_k}$$ as $G$ vector bundles.
 This comes from the fact that the orbits $\Gamma_i$
are $G$ stable, hence $G$ induces an action on each $X^{n_i}$.
Hence, the obstruction bundle splits as
$$E_{h_1,h_2} = \prod_{i=1}^k \left( H^1(\Sigma) \otimes
\Delta_i^*TX^{n_i} \right) ^G$$

We can simplify  the previous expression a bit further. Let $G_i$
be
 the subgroup of $S_{n_i}$ obtained from $G$ when its action is restricted to
the elements in $\Gamma_i$; then we have a surjective homomorphism
$$\lambda_i : G \rightarrow G_i$$ where the action of $G$ into
$\Delta_i^*TX^{n_i}$ factors through $G_i$. So we have

$$\left( H^1(\Sigma) \otimes \Delta_i^*TX^{n_i} \right) ^G \cong
\left( H^1(\Sigma)^{ker(\lambda_i)} \otimes \Delta_i^*TX^{n_i}
\right) ^{G_i}.$$

Now let $\Sigma_i:= \Sigma / \ker(\lambda_i)$, it is an orbifold
Riemann surface with a $G_i$ action so that  $\Sigma_i /G_i$
becomes an orbifold sphere with three marked points (the markings
are with respect to the generators of $G_i$: $h_1, h_2,
(h_1h_2)^{-1}$). So, in the same way as in the definition of the
obstruction bundle $E_{({\bf h})}$ we get that
$$E_{h_1,h_2}^i:=\left( H^1(\Sigma_i) \otimes \Delta_i^*TX^{n_i}
\right) ^{G_i}$$

The obstruction bundle splits as
$$E_{h_1,h_2}=\prod_{i=1}^k E_{h_1,h_2}^i$$

As the action of $G_i$ in $\Delta_i^*TX^{n_i}$ is independent on
the structure of $X$( moreover, it depends only in the
coordinates), hence

$$\Delta_i^*TX^{n_i} \cong TX \otimes \C^{n_i}$$
as $G_i$-vector bundles, where $TX$ is the tangent bundle over $X$
and $G_i \subset S_{n_i}$ acts on $\C^{n_i}$ in the natural way.
Then

$$E_{h_1,h_2}^i \cong TX \otimes (H^1(\Sigma) \otimes
\C^{n_i})^{G_i}$$

Defining $r(h_1,h_2)(i):=dim_\C (H^1(\Sigma) \otimes
\C^{n_i})^{G_i}$ it follows that the Euler class of
$E_{h_1,h_2}^i$ equals the Euler class of $X$ to some exponent
$c(E_{h_1,h_2}^i) = e(X)^{r(h_1,h_2)(i)}$. However, the underline
space is one copy of $X$ only. Hence,
$$c(E_{h_1,h_2}^i) = \left\{ \begin{array}{cc}
1 & \mbox{if }{r(h_1,h_2)(i)}=0 \\
e(X) & \mbox{if } {r(h_1,h_2)(i)}=1\\
0 & \mbox{if }{r(h_1,h_2)(i)} \geq 2
 \end{array} \right. \leqno(3.8)$$

    Therefore, we prove that
    \vskip 0.1in
    \noindent
    {\bf Theorem 3.2: }{\it
$$c(E_{h_1,h_2}) = \prod_{i=1}^k c(E^i_{h_1,h_2}),$$
where $c(E^i_{h_1,h_2})$ is given by previous formula (3.8).}
    \vskip 0.1in
    Theorem 3.2 matches the calculation of Lehn-Sorger for the
    cohomology of Hilbert scheme of points of $K3, T^4$ except a
    sign (see section 3.3).

\subsection{Conjecture for general crepant resolution}

Suppose that $\pi: Y\rightarrow X$ is one crepant resolution of Gorenstein  orbifold $X$.  Then, $\pi$ is
	a Mori contraction and the homology classes of rational curves $\pi$ contracted are generated by so called
	extremal rays. Let $A_1, \cdots,
    A_k$ be an integral basis of extremal rays. We call $\pi$ non-degenerate if  $A_1, \cdots,
    A_k$ are linearly independent. For example, the Hilbert-Chow
	map $\pi: M^{[n]}\rightarrow M^n/S_n$ satisfies this hypothesis.  Then,
    the homology class of any effective curve being contracted can be written as
    $A=\sum_i a_i A_i$ for $a_i\geq 0$. For each $A_i$, we assign a formal variable $q_i$. Then, $A$ corresponds to $q^{a_1}_1 \cdots q^{a_k}_k$.
    We define a 3-point function
    $$<\alpha, \beta, \gamma>_{qc}(q_1, \cdots q_k)=\sum_{a_1, \cdots, a_k}\Psi^X_{A}(\alpha, \beta, \gamma)q^{a_1}_1\cdots
    q^{a_k}_k ,\leqno(3.9)$$
    where $\Psi^X_{A}(\alpha, \beta, \gamma)$ is Gromov-Witten invariant and $qc$ stands for the quantum correction.
    We view $<\alpha, \beta, \gamma>_{qc}(q_1, \cdots, q_k)$
    as  analytic function of $q_1, \cdots q_k$ and  set $q_i=-1$ and let
    $$<\alpha, \beta,\gamma>_{qc}=<\alpha, \beta,\gamma>_{qc}(-1, \cdots, -1).\leqno(3.10)$$
    We define a quantum corrected triple intersection
    $$<\alpha, \beta, \gamma>_{\pi}=<\alpha, \beta,
    \gamma>+<\alpha, \beta, \gamma>_{qc},$$
    where $<\alpha, \beta, \gamma>=\int_X \alpha\cup\beta\cup\gamma$ is the ordinary triple intersection.
    Then we define the quantum corrected cup product $\alpha\cup_{\pi}
    \beta$ by the equation
    $$<\alpha\cup_{\pi} \beta, \gamma>=<\alpha, \beta,
    \gamma>_{\pi},$$
    for arbitrary $\gamma$. Another way to understand $\alpha\cup_{\pi}\beta$ is as following.
    Define a product $\alpha\star_{qc}\beta$ by the equation
    $$<\alpha\star_{qc}\beta, \gamma>=<\alpha, \beta,\gamma>_{qc}$$
    for arbitrary $\gamma \in H^*(Y, \C)$. Then, the quantum corrected product is the ordinary cup product
    corrected by $\alpha\star_{qc}\beta$. Namely,
    $$\alpha\cup_{\pi} \beta=\alpha\cup \beta +\alpha\star_{qc} \beta. \leqno(3.11) $$
    We denote the new quantum corrected cohomology ring as $H^*_{\pi}(Y, \C)$.

    \vskip 0.1in
    \noindent
    {\bf Cohomological Crepant Resolution Conjecture: }{\it Suppose that $\pi$ is non-degenerate and hence
	$H^*_{\pi}(Y, \C)$ is well-defined. Then,   $H^*_{\pi}(Y, \C)$
              is the ring isomorphic to orbifold cohomology ring $H^*_{orb}(X, \C)$.}
    \vskip 0.1in
    Recently, Li-Qin-Wang \cite{LQW3} proved a striking theorem that the
    cohomology ring of $M^{[n]}$ is universal in the sense that it
    depends only on homotopy type of $M$ and $K_X$. Combined with
    their result,  CCRC yields
    \vskip 0.1in
    \noindent
    {\bf Conjecture 3.3: }{\it $\star$ product depends only on $K_M$.}
    \vskip 0.1in
    It suggests an interesting way to calculate $\star$ product by
    first finding a universal formula (depending only on $K_M$)
    and calculating a special example such as $\P^2$ to determine
    the coefficient.

    Next, we formulate a closely related conjecture for $K$-equivalent manifolds.

    Suppose that $X, X'$ are $K$-equivalent and $\pi: X\rightarrow X'$ is the birational map. Again, exceptional rational curves makes sense.
     Suppose that $\pi$ is nondegenerate. Then, we go through the previous
    construction to define ring $H^*_{\pi}(X, \C)$.
    \vskip 0.1in
    \noindent
    {\bf Cohomological Minimal Model Conjecture: }{\it Suppose that $\pi, \pi^{-1}$ are nondegenerate. Then,
 $H^*_{\pi}(X, \C)$ is the ring isomorphic to $H^*_{\pi^{-1}}(X', \C)$ }
    \vskip 0.1in
    When $X, X'$ are the different crepant resolutions of the same orbifolds, Cohomological minimal model conjecture follows from Cohomological
    crepant resolution conjecture. However, it is well-known that most of K-equivalent manifolds are not crepant resolution of orbifolds.
    Cohomological minimal model conjecture can be generalized to orbifold provided that the quantum corrections are defined using orbifold
    Gromov-Witten invariants introduced by Chen-Ruan \cite{CR2}.
    \vskip 0.1in
    \noindent
    {\bf Remark 3.4: }{\it (1) The author does not know how to define quantum corrected product if $\pi$ is not
	nondegenerate. (2) All the conjectures in this section
    should be understood as the conjectures up to certain slight
    modifications (see next section).}
\vskip 0.1in
    \noindent
    {\bf Example 3.5: } Next, we use the work of Li-Qin \cite{LQ}
    to verify Cohomological Crepant Resolution Conjecture for
    $M^{[2]}$. To simplify the formula, we assume that $M$ is simply
    connected.

    It is easy to compute the orbifold cohomology $H^*_{orb}(X,
    \C)$ for $X=M^2/\Z_2$. The nontwisted sector
    can be identified with invariant cohomology of $M^2$.
    Let $h_i\in H^2(M, \C)$ be a basis and $H\in H^4(M, \C)$ be Poincare dual to
    a point. Then, the cohomology of the nontwisted
    sectors are generated by $1, 1\otimes h_i+h_i\otimes 1, 1\otimes
    H+H\otimes 1, h_i\otimes h_j+h_j\otimes h_i, h_i\otimes H+H\otimes h_i,
    H\otimes H.$ The twisted sector is diffeomorphic to
    $M$ with degree shifting number 1. We use $\tilde{1},
    \tilde{h_i}, \tilde{H}$ to denote the generators. They are of
    degrees $2,4,6$. By the definition, triple intersections
    $$<twisted sector, nontwisted sector, nontwisted sector>=0,$$
    $$<twisted sector, twisted sector, twisted sector>=0.$$
    Following is the table of nonzero triple intersections
    involving classes from the twisted sector
    $$<\tilde{1}, \tilde{1}, 1\otimes H+H\otimes 1>=1, <\tilde{1},
    \tilde{1}, h_i\otimes h_j+h_j\otimes h_i>=<h_i, h_j>, $$
    $$ <\tilde{1}, \tilde{h}_i, h_j\otimes 1+h_j\otimes 1>=<h_i, h_j>. \leqno(3.12)$$

    Next, we review the construction of $Y=M^{[2]}$.
    Let $\widetilde{M^2}$ be the blow-up of
    $M^2$ along the diagonal. Then, $\Z_2$ action
    extends to $\widetilde{M^2}$. Then, $Y=\widetilde{M^2}/\Z_2$.
    It is clear that we should map the classes from nontwisted sector to its pull-back $\pi: Y\rightarrow X$.
    We use the same notation to denote them. The exceptional divisor $E$ of Hilbert-Chow map
    $\pi: Y\rightarrow X$ is a $P^1$-bundle over $M$. Let
    $\bar{1}, \bar{h}_i, \bar{H}$ be the Poincare dual to $E$,
     $p^{-1}(PD(h_i))$ and fiber $[C]$, where $p: E\rightarrow M$ is the projection.

    Notes that $\bar{1}|_E=2\E$, where $\E$ is the tautological
    divisor of $P^1$-bundle $E\rightarrow M$. It is clear that
    $E=P(N_{\Delta(X)|X^2})$, where $\Delta(X)\subset X^2$ is the diagonal. Hence,
    $$<\bar{1}, \bar{1},
    \bar{h}_i>=4\E^2|_{p^{-1}(PD(h_i))}=4C_1(N_{\Delta(X)|X^2})\E|_{p^{-1}(PD(h_i))}=-4<C_1(X),
    h_i>. \leqno(3.13)$$
    $$<\bar{1}, \bar{1}, 1\otimes H+H\otimes 1>=2\E(1\otimes
    H+H\otimes 1)(E)=4\E(C)=-4.\leqno(3.14)$$
    $$<\bar{1}, \bar{1}, h_i\otimes h_j+h_j\otimes h_i>=-4<h_i,
    h_j>.\leqno(3.15)$$
    $$<\bar{1}, \bar{h}_i, 1\otimes h_j+h_j\otimes 1>=-4<h_i,
    h_j>.\leqno(3.16)$$
    Others are zero.

    The quantum corrections have been computed by Li-Qin
    \cite{LQ} (Proposition 3.021). The only nonzero terms are
    $$\begin{array}{lll}
    <\bar{1}, \bar{1}, \bar{h}_i>_{qc}(q)&=&\sum_{d=1}
    \bar{1}(d[C])^2\Psi^X_{d[C]}(\bar{h})q^d\\
    &=&\sum_{d=1} \frac{4d^2(2<K_X, h_i>)}{d^2} q^d\\
    &=&8<K_X, h_i>\frac{q}{1-q}
    \end{array}.\leqno(3.17)$$
    Hence,
    $$<\bar{1}, \bar{1}, \bar{h}_i>_{qc}=4<K_X, h_i>=4<C_1(X),
    h_i>\leqno(3.18)$$
    cancels $<\bar{1}, \bar{1},\bar{h}_i>$.

    It is clear that the map $\tilde{1}\rightarrow 2\bar{1},
    \tilde{h}\rightarrow 2\bar{h}, \tilde{H}\rightarrow \bar{H}$ is
    a ring isomorphism. $\Box$

    Next, we give two examples to verify  Cohomological Minimal Model Conjecture (CMMC).
    \vskip 0.1in
    \noindent
    {\bf Example 3.6: } The first example is the flop in dimension
    three. This case has been worked out in great detail by Li-Ruan
    \cite{LR}. For example, they proved a theorem that quantum cohomology rings are
    isomorphic under the change of the variable $q\rightarrow
    \frac{1}{q}$. Notes that if we set $q=-1$, $\frac{1}{q}=-1$. We set other quantum variables zero.
    Then, the quantum product becomes the quantum corrected product
    $\alpha\cup_{\pi}\beta$. Hence, CMMC follows from Li-Ruan's
    theorem. However, It should be pointed out that one can directly
    verify CMMC without using Li-Ruan's theorem. In fact, it is
    an much easier calculation.

    \vskip 0.1in
    \noindent
    {\bf Example 3.7: } There is a beautiful four dimensional
    birational
    transformation called Mukai transform as follows. Let
    $\P^2\subset X^4$ with $N_{\P^2|X^4}=T^* \P^2$. Then, one can
    blow up $\P^2$. The exceptional divisor of the blow up is  a hypersurface of $\P^2\times \P^2$
    with the bidegree $(1,1)$. Then, one can blow down in another
    direction to obtain $X'$. $X, X'$ are $K$-equivalent. In his
    Ph.D thesis \cite{Z}, Wanchuan Zhang showed that the quantum
    corrections $<\alpha, \beta, \gamma>_{qc}$ are trivial, and
    cohomologies of $X, X'$ are isomorphic.

    \subsection{Miscellaneous issues}

In the computation of orbifold cohomology of symmetry product
    and its relation to that of Hilbert scheme of points, there
    are two issues arisen. It was showed in the work of
    Lehn-Sorger, Fantechi-G\"{o}ttsche and Uribe that one has to
    add a sign in the definition of orbifold product in order to
    match that of Hilbert scheme of points. This sign was
    described as follows.

    Recall the definition of orbifold cup
    product
    $$\alpha\cup_{orb}\beta =\sum_{(h_1, h_2)\in T_2, h_i\in
    (g_i)} (\alpha\cup_{orb}\beta)_{(h_1,h_2)},\leqno(3.19)$$
    where $(\alpha\cup_{orb}\beta)_{(h_1,h_2)}\in
    H^*(X_{(h_1h_2)}, \C)$ is defined by the relation
    $$<(\alpha\cup_{orb}\beta)_{(h_1,h_2)}, \gamma>_{orb}
    =\int_{X_{(h_1,h_2)}}e^*_1\alpha\wedge e^*_2\beta\wedge e^*_3\gamma \wedge
    e_A(E_{(\g)}).\leqno(3.20)$$
    for $\gamma\in H^*_c(X_{((h_1h_2)^{-1})}, \C)$
    Then, we add a sign to each term.
    $$\alpha\cup_{orb}\beta =\sum_{(h_1, h_2)\in T_2, h_i\in
    (g_i)} (-1)^{\epsilon(h_1, h_2)}(\alpha\cup_{orb}\beta)_{(h_1,h_2)},\leqno(3.21)$$

    where
    $$\epsilon(h_1,
    h_2)=\frac{1}{2}(\iota(h_1)+\iota(h_2)-\iota(h_1h_2)).\leqno(3.22)$$
    Since
    $$\epsilon(h_1, h_2)+\epsilon(h_1h_2,
    h_3)=\frac{1}{2}(\iota(h_1)+\iota(h_2)+\iota(h_3)-\iota(h_1h_2h_3))=\epsilon(h_1,
    h_2h_3)+\epsilon(h_2, h_3),$$
    such a sign modification does not affect the associativity of
    orbifold cohomology.

	However, Qin-Wang \cite{QW} observed that the orbifold cohomology modified by such a sign
	is isomorphic to original orbifold cohomology over complex number by an explicit isomorphism
 	$$\alpha\rightarrow (-1)^{\frac{\iota(g)}{2}} \alpha$$
	for $\alpha\in H^*(X_{(g)}, \C)$. $\epsilon(h_1, h_2)$ is often an integer (for
	example symmetric product) while $\frac{\iota(g)}{2}$ is just a fraction. Hence, 
	$(-1)^{\frac{\iota(g)}{2}} $ is a complex number only.

    Another issue is the example of the crepant resolution of surface singularities
      $\C^2/\Gamma$. As Fantechi-G\"{o}ttsche \cite{FG2} pointed out, the Poincare paring of $H^2_{orb}(\C^2/\Gamma, \C)$
    is indefinite while the Poincare paring of its crepant resolution is negative definite. There is an
    easy way to fix this case (suggested to this author by Witten). We view the involution $I: H^*(X_{g^{-1}}, \C)\rightarrow H^*(X_g, \C)$
    as a "complex conjugation". Then, we define a "hermitian inner product"
    $$<<\alpha, \beta>>=<\alpha, I^*(\beta)>.\leqno(3.23)$$
    If we use this "hermitian" inner product, the intersection paring is positive definite again.
     The above process has its conformal theory origin (see
     \cite{NW}). It is attempting to perform this modification on
     $H^*(X_{(g)}, \C)\oplus H^*(X_{(g^{-1})}, \C)$ whenever
     $\iota_{(g)}=\iota_{(g^{-1})}$.
  The author does not know if it will affect  the associativity of
      orbifold cohomology.

      In a different direction, it is useful to modify
    orbifold cup product by a "complex" sign for following reason. In ordinary
    cohomology, cup product is supercommutative
    $$\alpha\cup
    \beta=(-1)^{deg(\alpha)deg(\beta)}\beta\cup \alpha.$$
    For the orbifold cup product, we only have such
    supercommutativity when degree shifting number $\iota(g)$ is
    integer. In general case, degree shifting number is a rational
    number and we are talking about $(-1)^a$ for a fraction $a$.
    At the first, it may sounds odd since   "super" means odd/even property. However, there were
    other reasons to believe that it is important to consider a
    "complex sign" $(-1)^a$ even for a fraction $a$. Here, we
    interpret $(-1)^a$ as a complex number $e^{\pi i a}$. The author does not
    know how to modify the orbifold product to achieve
    such a "complex supercommutativity".

      \section{Twisted orbifold K-theory}

      An important aspect of the stringy orbifold is the orbifold
      K-theory and its twisted version. Unfortunately, the twisted
      version developed by Adem-Ruan \cite{AR} is less than satisfactory. As
      showed in \cite{R3}, one can perform the twisting of
      the orbifold cohomology for any inner local system. However, the
      twisted orbifold K-theory has only been constructed for
      the discrete torsion. The decomposition theorem by Chern homomorphism in the twisted case was only proved
      for the global quotient although the author expects it to be
      true in general. As Adem and the author was developing our
      twisted orbifold K-theory, they noted a parallel development
      of the twisted K-theory on smooth manifolds by Witten and others \cite{W}
      for the purpose of
      describing the D-brane charge in physics. They were wondering if
      one can unify these two twisting. This was done beautifully
      by Lupercio-Uribe \cite{LU1} using the notion of gerbes,
      which is an interesting concept of its own. Their construction will be presented in this section.

      However, the
      corresponding decomposition theorem or the existence of corresponding Chern homomorphism
      is far less obvious. For
      example, we do not yet know the analogous of the orbifold cohomology
      for Witten's twisted K-theory. Since Lupercio-Uribe's
      construction works over more general orbit space or the Artin
      stack,
      it suggests that we should have a
      version of the orbifold cohomology for a general orbit space and
      its corresponding Chern homomorphism. There is a striking result recently by
      Freed-Hopkings-Teleman \cite{FR} which shows that the twisted equivariant K-theory of Lie group
      $G$ acting on itself by the conjugation is isomorphic to the Verlinde
      algebra. However, their calculation did not explain the
      reason of the appearance of Verlinde algebra, i.e, establishing the explicit isomorphism. I am sure that the
      existence of the corresponding "orbifold cohomology" and Chern
      homomorphism should clarify this issue. However, I would
      like to warn the reader that the straight forward
      generalization of Adem-Ruan's decomposition theorem does not work. In fact, there
      was an old article \cite{H} to show that such a Chern
      homomorphism could not exist if an action has $S^1$-isotropy
      subgroup.

      Another topic we will survey in this section is in different
      nature. Wang constructed the twisted K-theory using a Clifford
      algebra (Ferminon) \cite{W}. The generating function of its Euler
      characteristic has a better shape. It is an interesting
      construction. However, this author does not know how to
      corporate it into other type of twisting.

      \subsection{Gerbe and twisted K-theory}

        Let $\RR \stackrel{\rightarrow}{\rightarrow} \UU$ be groupoid associated to an orbifold $X$.
    \vskip 0.1in
    \noindent
    {\bf Definition 4.1: }{\it
A gerbe over an orbifold $\RR \stackrel{\rightarrow}{\rightarrow}
\UU$, is a complex line bundle $\LL$ over $\RR$ satisfying the
following conditions
\begin{itemize}
\item $i^* \LL \cong \LL^{-1}$
\item $\pi_1^* \LL \otimes \pi_2^* \LL \otimes m^*i^* \LL \stackrel{\theta}{\cong} 1$
\item $\theta : \RR_t\times_s \RR \to U(1)$ is a 2-cocycle
\end{itemize}
where $\pi_1,\pi_2 : \RR_t\times_s \RR \to \RR$ are  the
projections on the first and the second coordinates, and $\theta$
is a trivialization of the line bundle. }
    \vskip 0.1in
    When the groupoid $\G$ is given by the action of a finite group on
    a point $G\times \star\stackrel{\rightarrow}{\rightarrow}
    \star$, a gerbe over $\G$ is the same as an central extension
    $$1\rightarrow \overline{U(1)}\rightarrow \tilde{\G}\rightarrow
    \G\rightarrow 1.$$
     Gerbes over a discrete group $G$ are in 1-1
    correspondance with the set of two cycles $Z(G, U(1))$.

    The set of gerbes forms a group $Gb(\RR\stackrel{\rightarrow}{\rightarrow} \U)$.
    Gerbes has a characteristic class $<\L>$ given by $\theta:
    \RR_t\times_s \RR\rightarrow U(1)$. $<\L>$ is a degree three
    cohomology in so called classifying space of groupoid.
    Moreover, the group $Gb(\RR\stackrel{\rightarrow}{\rightarrow}
    \UU)$is independent of the Morita class of $\RR\stackrel{\rightarrow}{\rightarrow}
    \UU$.

    \vskip 0.1in
    \noindent
    {\bf Example 4.2: }
   Consider an inclusion of (compact Lie) groups $K\subset G$ and consider the groupoid
   $\GG$ given by the action of $G$ in $G/K$,
   $$ G/K \times G \stackrel{\rightarrow}{\rightarrow} G/K $$

   Observe that the stabilizer of $[1]$ is $K$ and therefore we have that the
   following groupoid
   $$ [1]\times K \stackrel{\rightarrow}{\rightarrow} [1] $$
   is Morita equivalent to the one above.

   From this we obtain
   $$ Gb(\GG) \cong H^3(K,\Z) $$

 For a smooth manifold $X$ we have that
$$ Gb(X) = [ X , BB\C^* ] $$
where $BB\C^* = B\P U(\HH)$ for a Hilbert space $\HH$.

Let us write $\overline{\P U(\HH)}$ to denote the groupoid $\star
\times \P U(\HH) \to \star$. We have the following
    \vskip 0.1in
    \noindent
    {\bf Proposition 4.3: }{\it
   For an orbifold $X$ given by a groupoid $\Xx$ we have
   $$ Gb(\Xx) = [\Xx , \overline{\P U(\HH)}] $$
   where $[\Xx , \overline{\P U(\HH)}]$ represents the Morita equivalence classes of morphisms
   from $\Xx$ to $\overline{\P U(\HH)}$}
   \vskip 0.1in
Just as in the case of a gerbe over a smooth manifold, we can do
differential geometry on gerbes over an orbifold groupoid
$\Xx=(\RR\stackrel{\rightarrow}{\rightarrow} \UU)$. Let us define
a connection over a gerbe in this context.
    \vskip 0.1in
    \noindent
    {\bf Definition 4.4: }{\it
   A connection $(g,A,F,G)$ over a gerbe  consists of a complex valued
   0-form $g\in\Omega^0 (\RR_t\times_s \RR)$, a 1-form $A\in \Omega^1(\RR)$, a 2-form $F\in\Omega^2(\UU)$ and a
   3-form $G\in\Omega^3(\UU)$ satisfying
   \begin{itemize}
      \item $G=dF$,
      \item $t^* F - s^* F = dA$ and
      \item $\pi_1^* A + \pi_2^* A +  m^* i^* A = -\sqrt{-1} g^{-1} dg$
   \end{itemize}
   The 3-form $G$ is called the $curvature$ of the connection. A connection is called
   $flat$ if its curvature $G$ vanishes.}
   \vskip 0.1in

The 3-curvature $\frac{1}{2\pi\sqrt{-1}}G$ represents the integer
characteristic class of the gerbe in cohomology with real
coefficients, this is the Chern-Weil theory for a gerbe over an
orbifold. Recall that for a line bundle over a manifold $X$ with a
connection, the holonomy map can be considered as a map  from loop
space $LX$ to $U(1)$. For gerbes over $X$ with a connection, there
is a notion of holonomy. But it produces a flat line bundle over
$LX$.   If we are given a gerbe $L$ with a flat connection over a
groupoid $\Xx$, using the holonomy we construct a "flat line
bundle" $\Lambda$ over the loop groupoid $\LL\Xx$ and hence obtain
a flat line bundle over the fixed point set $\wedge
\Xx=\LL\Xx^{S^1}$ of the natural circle action.

    \vskip 0.1in
    \noindent
    {\bf Theorem 4.5: }{\it The inertia orbifold
 $\widetilde{\Sigma_1 X}$ defined in \cite{CR1}
 is represented by the groupoid $\wedge \Xx$.  The holonomy line
bundle $\Lambda$ over $\wedge \Xx$ is an inner local system as
defined in \cite{R3}.}
    \vskip 0.1in
    However, we do not know if every inner local system can be obtained in this way although
    it can recover the inner local system induced by discrete torsion (see the comments after
    Proposition 4.8).

    Next, we use gerbes to construct twisted K-theory. We first
    present the construction of "twisted" vector bundle, a
    generalization of Adem-Ruan twisted orbifold vector bundle.
    But this construction only yields nontrivial K-theory when
    $<\L>$ is a torsion class. For non-torsion class, we present a
    construction generalizing both twisted vector bundle and
    twisted K-theory in the smooth case.
    \vskip 0.1in
    \noindent
    {\bf Definition 4.6: }{\it
An $n$-dimensional $\LL$-twisted bundle over $\RR
\stackrel{\rightarrow}{\rightarrow} \UU$ is a vector bundle $E
\rightarrow \UU$
   together with a given isomorphism $$\LL \otimes t^* E \cong s^* E
   $$.}
    \vskip 0.1in

Notice that we then have a canonical isomorphism
$$ m^* \LL \otimes \pi_2^*t^* E \cong \pi_1^* \LL \otimes \pi_2 ^* \LL
\otimes \pi_2^* t^* E \cong \pi_1^* \LL \otimes \pi_2^*(\LL
\otimes t^* E) \cong \pi_1 ^* \LL \otimes \pi_2^* s^* E$$

    We can define the corresponding Whiteney sum of $\LL$-twisted
    bundles similarly.
    \vskip 0.1in
    \noindent
    {\bf Definition 4.7: }{\it
The Grothendieck group generated by the isomorphism classes of
$\LL$ twisted bundles over the orbifold $X$ together with the
addition operation just defined is called the $\LL$ twisted
$K$-theory of $X$ and is denoted by $^\LL K_{gpd}(X)$.}
    \vskip 0.1in

Moerdijk and Pronk \cite{MP1} proved that the isomorphism classes
of orbifolds are in 1-1 correspondance with the classes of
\'etale, proper groupoids up to Morita equivalence.
 As a direct consequence of the definitions,
 $^\LL K_{gpd}(X)$ is independent of the groupoid
that is associated to $X$.

    Using the group structure of $Gb(\RR \stackrel{\rightarrow}{\rightarrow} \UU)$
we can define a product between bundles twisted by different
gerbes, so for $\LL_1$ and $\LL_2$ gerbes over $X$
$$^{\LL_1} K_{gpd}(X) \otimes  ^{\LL_2} K_{gpd}(X) \rightarrow ^{\LL_1\otimes\LL_2} K_{gpd}(X)$$
and we can define the total twisted orbifold $K$-theory of $X$ as
$$T K_{gpd}(X) = \bigoplus_{\LL \in Gb(\RR \stackrel{\rightarrow}{\rightarrow} \UU)} ^\LL K_{gpd}(X)$$
This has a ring structure due to the following proposition.

    \vskip 0.1in
    \noindent
    {\bf Proposition 4.8: }{\it
   The twisted groups $^\LL K_{gpd}(\GG)$ satisfy the following properties:
   \begin{itemize}
   \item[1.] If $\langle \LL \rangle =0$ then $^{\LL} K_{gpd}(\GG) = K_{gpd}(\GG)$.
   \item[2.] $^\LL K_{gpd}(\GG)$ is a module over $^\LL K_{gpd}(\GG)$
   \item[3.] If $\LL_1$ and $\LL_2$ are two gerbes over $\GG$ then there is a
   homomorphism
   $$^{\LL_1} K_{gpd}(\GG) \otimes  ^{\LL_2} K_{gpd}(\GG) \rightarrow ^{\LL_1\otimes\LL_2} K_{gpd}(\GG)$$
   \item[4.] If $\psi\colon \GG_1 \longrightarrow \GG_2$ is a groupoid homomorphism
   then there is an induced homomorphism
   $$^\LL K_{gpd}(\GG_2) \longrightarrow ^{\psi^* \LL} K_{gpd}(\GG_1) $$
   \end{itemize}.}
   \vskip 0.1in
   In  the case when $Y$ is the orbifold universal cover of $X$ with orbifold
fundamental group $\pi_1^{orb}(X) =H$,
 we  can take a discrete torsion $\alpha \in H^2(H,U(1))$ and define the twisted $K$-theory
of $X$ as in \cite{AR}. The groupoid $\RR_Y \times H
\stackrel{\rightarrow}{\rightarrow} \UU_Y$ represents the orbifold
$X$.  We want to construct a gerbe $\LL$ over $\RR_Y \times H
\stackrel{\rightarrow}{\rightarrow} \UU_Y$ so that the twisted
$^\LL K_{gpd}(X)$ is the same as the twisted $^\alpha K_{orb}(X)$
as defined by Adem-Ruan \cite{AD}.

The discrete torsion $\alpha$ defines a central extension of $H$
$$ 1 \rightarrow U(1) \rightarrow \widetilde{H} \rightarrow H \rightarrow 1$$
and doing the cartesian product with $\RR_Y$ we get a line bundle

\begin{displaymath}
\begin{array}{ccccc}
\overline{U(1)} & & \rightarrow &\LL_\alpha = & \RR_Y \times \widetilde{H} \\
 & & & & \downarrow \\
 & & & & \RR_Y \times H
\end{array}
\end{displaymath}
which ia a gerbe over $\RR_Y \times H
\stackrel{\rightarrow}{\rightarrow} \UU_Y$ .
    \vskip 0.1in
    \noindent
    {\bf Theorem 4.9: }{\it
$^{\LL_{\alpha}} K_{gpd}(X) \cong ^\alpha K_{orb}(X)$.}
    \vskip 0.1in
     We should point out here that the theory so far described is
essentially empty whenever the characteristic class $\langle \LL
\rangle$ is a non-torsion element in $H^3(M,\Z)$. The following is
true.
    \vskip 0.1in
    \noindent
    {\bf Proposition 4.10: }{\it
If there is an $n$-dimensional $\LL$-twisted bundle over the
groupoid $\GG$ then $\langle\LL\rangle ^n = 1 $.}
    \vskip 0.1in
    We need to consider a more general definition when
the class $\langle \LL \rangle$ is a non-torsion class. With this
in mind we would like to have a group model for the space $\FF$ of
Fredholm operators. One possible candidate is the following.
    \vskip 0.1in
    \noindent
    {\bf Definition 4.11: }{\it For a given Hilbert space $\HH$
by a {\em polarization of $\HH$} we mean a decomposition
$$\HH =\HH_+ \oplus \HH_- $$ where $\HH_+$ is a complete infinite
dimensional subspace of $\HH$ and $\HH_-$ is its orthogonal
complement.

We define the group $GL_{res}$ to be the subgroup of automorphism
of $\HH$ consisting of operators $A$ that when written with
respect to the polarization $\HH_+ \oplus \HH_-$ look like
$$A = \left(
      \begin{array}{cc}
      a & b \\
      c & d
      \end{array}
      \right)
$$
where $a\colon \HH_+ \rightarrow \HH_+$ and $d \colon \HH_-
\rightarrow \HH_-$ are Fredholm operators, and $b \colon \HH_-
\rightarrow \HH_+$ and $c \colon \HH_+ \rightarrow \HH_-$ are
Hilbert-Schmidt operators.}
    \vskip 0.1in

Notice that $\P U{\HH_+}$ acts by conjugation on $GL_{res}$ for
the map that sends $g\in \P U{\HH_+}$ to
  $$ \iota(g) = \left(
  \begin{array}{cc}
      g & 0 \\
      0 & 1
  \end{array}
  \right)
  $$

We have the following fact.
    \vskip 0.1in
    \noindent
    {\bf Proposition 4.12: }{\it The map $ GL_{res} \to \FF \colon A \rightarrow a $
is a homotopy equivalence.}
    \vskip 0.1in

 Consider  a gerbe $\LL$ with characteristic class $\alpha$ as a
map $ X \rightarrow BB U(1) = B \P U(\HH_+) $, namely a Hilbert
projective bundle $Z_\alpha(M) \rightarrow M$. Then we form a
$GL_{res}$-adjoint bundle  over $X$ by defining
$$\FF_\alpha(M) = Z_\alpha(M) \times_{ \P U{\HH}} GL_{res}.$$

    \vskip 0.1in
    \noindent
    {\bf Definition 4.13: }{\it We define the twisted K-theory
    $$^\LL K_{gpd}(X)$$ as the homotopy class of the sections of
    $\FF_{\alpha}(X)$.}
    \vskip 0.1in

This definition works for a gerbe whose class is non-torsion and
has the obvious naturality conditions. In particular it becomes
Witten's twisting if the groupoid represents a smooth manifold.

      \subsection{Twisting by fermion}
          It is known that $H^2(S_n, U(1))=\Z_2$ for $n\geq 4$ and
          there is a twisted orbifold K-theory for symmetry
          product. The generating function of its Euler
          characteristic has been computed by Dijkgraaf \cite{D},
          Wang \cite{W} and Uribe \cite{U1}. Unlike the generating
          function of Euler characteristic of ordinary orbifold
          K-theory, this twisted version does not have modularity.
          To remedy the situation, Wang \cite{W} developed a super
          version of (equivariant) twisted K-theory for global
          quotient whose Euler characteristic again has nice
          property.  It is a very interesting construction. Let me
          sketch his construction.

          Let $\tilde{G}$ be a finite group and  $d:
          \tilde{G}\rightarrow \Z_2$ be a group epimorphism
          understood as parity function. An element $a$ of $\tilde{G}$
          is called odd/even if $d(a)=1/0$. $(\tilde{G}, d)$ is
          called a {\em finite supergroup}. Recall that spin
          group is double cover of $SO(n)$. We also require
          $\tilde{G}$ is double cover of $G$ and the distinguished
          central element $\theta$ of order 2 is even. Next, we
          consider its representation theory. Again, a
          representation of a finite supergroup is the same as a
          module of its group superalgebra. Moreover, we only
          consider supermodule, i.e., $\Z_2$-graded.  Given two
          supermodules $M=M_0+M_1$ and $N=N_0+N_1$ over a
          superalgebra $A=A_0+A_1$, the linear map $f:
          M\rightarrow N$ between two $A$-supermodules is a
          homomorphism of degree $i$ if $f(M_j)\subset M_{i+j}$
          and for any homogeneous element $a\in A$ and any
          homogeneous vector $m\in M$ we have
          $$f(am)=(-1)^{d(f)d(a)}a f(m).$$
          The degree $0/1$ part of a superspace is referred to as
          the even/odd part. We denote
          $$Hom_A(M, N)=Hom_A(M, N)_0\oplus Hom_A(M, N)_1,$$
          where $Hom_A(M, N)_i$ consist of $A$-homomorphisms of
          degree $i$ from $M$ to $N$. The notions of submodules,
          tensor product, and irreducibility for supermodules are
          defined similarly. Given a finite supergroup
          $\tilde{G}$, a $\tilde{G}$-supermodule $V$ is called
          {\em spin} if the central element $\theta$ acts as $-1$,
          a reministic of spin representation. We will only
          consider spin supermodule. A spin supermodule can also
          be thought as a $\alpha$-twisted projective module for a
          nontrivial 2-cocycle $\alpha: G\times G\rightarrow
          \Z_2$.

          There are two types of complex simple superalgebras
          $M(r|s)$ and $Q(n)$. $M(r|s)$ is the superalgebra
          consisting of the linear transformations of the
          superspace $\C^{r|s}=C^r+\C^s$. The superalgebra $Q(n)$
          is the graded subalgebra of $M(n|n)$ consisting of
          matrices of the form
          $$\left(\begin{array}{cc}
                    A&B\\
                    B&A
                    \end{array}\right)$$
          According to the classification of simple superalgebras
          above, the irreducibe supermodules of a finite
          supergroup are divided into two types, type $M$ and type
          $Q$. We note that the endomorphism algebra of an
          irreducible supermodule $V$ is isomorphic to $\C$ if $V$
          is of type $M$ and isomorphic to the complex Clifford
          algebra $C_1$ in one variable if $V$ is of type $Q$.

          For any conjugacy class $C$ of $G$, $\theta^{-1}(C)$ is
          either a conjugacy class of $\tilde{G}$ or it splits
          into two conjugacy classes of $\tilde{G}$. $C$ (and any element
          $g\in C$) is called  split or non-split. It is easy to
          check that $g$ is split iff the twisted character
          $L^{\alpha}_g=\alpha(g, .)\alpha(g,.)^{-1}$ is a trivial
          character of the centralizer $C(g)$.

          Now, we are ready to globalize above construction.
          Suppose that a finite group $G$ acts on a smooth
          manifold $X$ so that $X/G$ is a global quotient
          orbifold. Let $\tilde{G}$ be corresponding supergroup. A
          spin $\tilde{G}$-vector superbundle $E$ has property
          $E=E_0+E_1$ (often denoted by $E_0|E_1$) over $X$ such
          that $\tilde{G}$ acts in a $\Z_2$-graded manner and
          $\theta$ acts as $-1$.  Since $\tilde{G}$ always
          contains odd elements which interchange $E_0, E_1$.
          Hence, $rank E_0=rank E_1$. The direct sum operation
          extends to spin $\tilde{G}$-bundle and we can define its
          K-theory $K^{-}_{\tilde{G}}(X)$. For this super-twisted
          K-theory, Wang proved a decomposition theorem similar to
          that of Adem-Ruan \cite{AR}.
          \vskip 0.1in
          \noindent
          {\bf Theorem 4.14: }{\it Let $X, \tilde{G}$ be define as
          above. We have a natural $\Z_2$-graded isomorphism
          $$K^{-}_{\tilde{G}}(X)\otimes \C\cong \oplus_{[g]}
          (K(X^g)\otimes L^{\alpha}_g)^{C(g)},\leqno(4.1)$$
          where the summation runs over the even conjugacy classes
          in $G$.}
          \vskip 0.1in
          If we drop super condition, we simply obtain twisted
          K-theory $^\alpha K_{G}(X)$ in the sense of Adem-Ruan.
          The corresponding decomposition theorem is a summation
          over all the conjugacy classes.
          \vskip 0.1in
          \noindent
          {\bf Example 4.15: }{\it Recall that the symmetry group $S_n$
          has nontrivial discrete torsion given by the fact that
          spin group is a double cover of $SO(n)$ and $S_n$ can be
          embedded into $SO(n)$ as a subgroup. This nontrivial
          discrete torsion can be described by exact sequence
          $$1\rightarrow \Z_2\rightarrow \tilde{S}_n\rightarrow
          S_n,$$
          where $\tilde{S}_n$ is generated by  $\theta$ and $t_i, i=1,
          \cdots, n-1$ and subject to the relations
          $$\theta^2=1, t^2_i=(t_it_{i+1})^3=\theta, t_it_j=\theta
          t_jt_i (i>j+1), zt_i=t_iz.$$
          The group $\tilde{S}_n$ carries a natural $\Z_2$ grading
          by letting $t_i$'s be odd and $\theta$ be even. Hence,
          $\tilde{S}_n$ is a finite supergroup. Wang studied
          $K^{-}_{\tilde{S}_n}(X^n)$ in detail and shows that
          $\oplus^{\infty}_{n=0}K^{-}_{\tilde{S}_n}(X^n)$ has many
          beautiful algebraic structure such as  Hopf algebra,
          admitting an action of the twisted Heisenberg algebra.
          We should emphasis that the ordinary twisted K-theory by
          discrete torsion does not have these nice algebraic
          structure.}

\section{Spin Orbifold Quantum Cohomology}

    One of the attractions of the stringy orbifold is its unique orbifold
    feature which will not exist on smooth manifolds. In the
    on-going project of this author with T. Jarvis, such a unique
    orbifold feature is added to its quantum cohomology theory.
    This is the incorporation of the theory of spin curves to
    the orbifold quantum cohomology. It is even more remarkable that
    the theory of spin curves should be naturally an orbifold theory and
    was discovered ten years  before its time. Now, it
    finally comes to its rightful home. This author should mention
    that an orbifoldish construction of the moduli space of spin
    curves has already been obtained by D. Abramovich and T.
    Jarvis \cite{AJ}. Here, we go one more step to incorporate it into the
    Gromov-Witten theory.

    The history of spin curves goes
    back ten years ago when Witten formulated a famous conjecture to
    relate  intersection numbers on the moduli space of stable Riemann
    surfaces to the KdV-hierarchy. Witten's conjecture was proved by
    Kontsevich \cite{K}, and motivated another well-known conjecture
    by Eguchi-Hori-Xiong and Katz in the quantum cohomology that the
    generating function of Gromov-Witten invariants satisfies a set of
    equations which form a Virasoro algebra. This conjecture is
    commonly known as the \emph{Virasoro conjecture}. Witten's
    conjecture was based on a matrix model description of so-called
    2D-gravity. Around the same time, Witten also formulated a closely
    related but less known conjecture for so-called \emph{2D-gravity
    coupled with topological matter}. This conjecture is equivalent to
    saying that the generating functions of intersection numbers on
    the moduli space of higher spin curves satisfies constraints
    forming a $W$-algebra---a larger algebra than the Virasoro
    algebra. Since then, his conjecture has been given a rigorous
    mathematical formulation by Jarvis and  his
    collaborators \cite{J1} \cite{J2} \cite{JKV} \cite{JKV2}.

    An incorporation of spin curves into quantum cohomology yields an elegant
    theory of spin Gromov-Witten theory which is again unique for the orbifold.
    Furthermore, its generating function should satisfy a set of differential
    equations forming a W-algebra. This is again an orbifoldish generalization of
    the Virasoro conjecture. Recall that the Gromov-Witten theory corresponds to
    the topological sigma model coupled with 2-D gravity in physics. In
    \cite{CR2}, it was showed that the sigma model on orbifold is
    not the ordinary map but good map which is the basis for orbifold
    Gromov-Witten invariants. It is amazing that the counterpart
    of 2-D gravity in the orbifold theory is not 2-D gravity but 2-D
    gravity coupled with the topological matter. This realization
    should has far reaching consequence to other area of interest
    such as mirror symmetry. For example, what is the B-model
    interpretation Spin orbifold quantum cohomology? A better
    understanding of this question should enrich the mirror symmetry
    as well.

    To illustrate the idea of the present paper, recall that the
    quantum cohomology can be thought of as a theory to study the
    Cauchy-Riemann equation:
    $$\bar{\partial}f=0.$$
    for a map $f: \Sigma\rightarrow X$, where $\Sigma$ is a Riemann
    surface and $X$ is a symplectic manifold. The orbifold quantum
    cohomology studies the same equation, where $\Sigma$ and $f$ are
    interpreted appropriately in an orbifold category. The spin
    orbifold quantum cohomology studies simultaneous solutions of the
    Cauchy-Riemann equation and another equation we call the
    \emph{Spin equation}.  The two equations are
    $$\bar{\partial}f=0,     \quad \bar{\partial}s+\bar{s}^{r-1}=0,$$
     where $s\in
    \Omega^0(L)$ for a line bundle $L$ which is an $r$th-root of the
    log canonical bundle (see Definition 5.7) of an orbifold
    Riemann surface. When $\Sigma$ is smooth and $r=2$, $L$
    corresponds a classical spin structure of $\Sigma$.  Hence, for
    arbitrary $r$, $L$ can be thought as a generalized spin structure.

    An easy computation (Lemma 5.1) shows that the second
    equation has only zero as its solution. So its solution (moduli)
    space is closely related to that of the Cauchy-Riemann equation.
    However, the spin equation is non-transverse in general; hence, it
    gives rather different invariants.

	Our spin orbifold quantum cohomology yields a new cohomology of
	orbifold, which we call r-spin orbifold cohomology or spin orbifold 
	cohomology. We do not yet know how useful spin orbifold cohomology is.

\subsection{Spin orbifold cohomology}
    One of remarkable consequences of ordinary orbifold quantum cohomology is
    the existence of a new cohomology ring of orbifold as the classical
    part of orbifold quantum cohomology. Our r-spin quantum cohomology also
    yields a new cohomology theory of orbifold, which is the main topic of
    another paper. In this section, we outline some of main ingredients.

    Let $X$ be an almost complex orbifold and $\B(\Z/r)$ be a point with the group
    action $\Z/r$. Let
    $X_{(g, l)}=X_{(g)}\times \B(\Z/r)_{(l)}$ for $0\leq l\leq r-1$. 

    \vskip 0.1in
	\noindent
	{\bf Definition 5.1.1: }{\it 
    We call  $X_{(g,l)}$ {\em Neveu-Schwarz}-sector
    for $l\neq 0$ and $X_{(g,0)}$ {\em Ramond sector}.}
    
    \vskip 0.1in
    \noindent
	{\bf Definition 5.1.2: }{\it 
    We define the $r$-spin orbifold cohomology $H^*_{r, orb}(X, \Q)$ of $X$ to be the
    vector space
    $$H^*_{r, orb}(X, \Q):=\oplus_{(g), l>0} H^*(X_{(g, l)}, \Q).$$}
	\vskip 0.1in
    
    As in the "usual" (without $r$-spin structure) case, there is a shift in the grading of the
    cohomology, and a Poincare paring. Let
    $$\iota(g, l):=\iota(g)+\frac{l-1}{r},$$
    where $\iota(g)$ is the usual degree shifting number  of $g$.
    The $r$-spin orbifold degree $deg_{r, orb}(\alpha)$ of $\alpha\in H^j(X_{(g,l)})$ is defined to be
    $$deg_{r, orb}(\alpha)=deg(\alpha)+2\iota(g,l).$$
    So
    $$H^d_{r, orb}(X, \Q)=\oplus_{(g), l>0}H^{d-2\iota(g)}(X_{(g, l)}, \Q).$$

   \vskip 0.1in
	\noindent
	{\bf Definition 5.1.3: }{\it 
    Let $I: X_{(g,l)}\rightarrow X_{(g^{-1}, r-l)}$ be the natural diffeormorphism. For all $0\leq d\leq
    2n+\frac{2(r-2)}{r}$ let
    $$<>_{r, orb}: H^{d}_{r, orb}(X, \Q)\otimes H^{2n+\frac{2(r-2)}{r}-d}_{r, orb}(X, \Q)\rightarrow \Q$$
    be defined as the direct sum of
    $$<>_{r, orb}^{(g,l)}: H^{d}(X_{(g,l)}, \Q)[-2\iota(g, l)]\otimes H^{2n+\frac{2(r-2)}{r}-d}
    (X_{(g^{-1},r-l)}, \Q)[-2\iota(g^{-1}, r-l)]\rightarrow \Q,$$
    where
          $$<\alpha, \beta>^{(g,l)}_{r, orb}=\int_{X_{(g, l)}}\alpha\wedge I^*(\beta).$$}
    \vskip 0.1in
	\noindent
	{\bf Remark 5.1.4: }{\it 
    The top degree of $H^*_{r, orb}(X, \Q)$ is $2n+\frac{2(r-2)}{r}$. Hence, it should be thought to
    have dimension $2n+\frac{2(r-2)}{r}$. In another words, $\B(\Z/r)$ contributes a dimension
    $\frac{2(r-2)}{r}$.}
	\vskip 0.1in
	\noindent
	{\bf Proposition 5.1.5: }{\it   $<>_{r, orb}$ is non-degenerate.}
	\vskip 0.1in
    
    The proof is essentially identical to the usual (non-spin) orbifold case.

    $H^*_{r, orb}(X, \Q)$ has a ring structure as the degree zero part of r-spin orbifold quantum
    cohomology we will define in the later section. 
	\vskip 0.1in
	\noindent
	{\bf Example 5.1.6: }{\it 
    Let $X$ be a point.  However, r-spin cohomology
    $H^*_{r, orb}(X, \Q)$ can be described as follows. Let $e_u$ for $u \in \{1,\dots,r-1\}$
    be the constant function $1$ on sector $X_{(1,u)}$. The $r$-spin cohomology is generated by
    $e_u$ with product  $e_u * e_v
    = e_{u+v-1}$ if $u+v<r+1$ and $e_u* e_v=0$ otherwise.}

    \subsection{Spin equation}

    Recall that a pseudo-holomorphic map $f: \Sigma\rightarrow X$ can
    be thought as a solution of the Cauchy-Riemann equation
    $\bar{\partial}f=0$. The solution space of $\bar{\partial}f=0$ is
    the moduli space of pseudo-holomorphic maps. The moduli space of
    stable maps is a ``nice" compactification of the moduli space of
    pseudo-holomorphic maps. A \emph{spin pseudo-holomorphic map} is a
    solution of two  nonlinear equations $$\bar{\partial}f=0, \quad
    \bar{\partial}s+\bar{s}^{r-1}=0.$$ The first equation is the usual
    one for pseudo-holomorphic maps, and the second equation uses a
    generalized spin structure on the Riemann surface as follows.
    Suppose that $\Sigma$ is a marked Riemann surface with marked points
    $z_1, \cdots, z_k$, and
    $\pi:\tilde{\Sigma}\rightarrow \Sigma$ is an orbifold Riemann
    surface with additional orbifold structure at marked point. If $L$ is an orbifold
    line bundle on $\tilde{\Sigma}$, it can be uniquely lifted to an orbifold
    line bundle over $\tilde{\Sigma}$. We shall use the same $L$ to
    denote its lifting.  Let    $K$ is the canonical bundle of $\Sigma$.

    \vskip 0.1in
	\noindent
	{\bf Definition 5.2.1: }{\it    Let $$K_{
    log} := K \otimes  \O(z_1) \otimes \dots \otimes
    \O(z_k)$$ be the \emph{log-canonical bundle}. $K_{log}$ can be thought as
    canonical bundle of punctured Riemann surface $\Sigma-\{z_1, \cdots, z_k\}$. Suppose that $L$ is an
    orbifold line bundle on $\tilde{\Sigma}$ with an isomorphism $\phi: L^r\rightarrow
    K_{log}$, where $K_{log}$ is identified with its pull-back on $\tilde{\Sigma}$.  The pair $(L, \phi)$ is called an $r$th-root, or a
    \emph{generalized spin structure}.}
    \vskip 0.1in
	\noindent
	{\bf Definition 5.2.2: }{\it 
    Suppose that the chart of $\tilde{\Sigma}$ at an orbifold point $z_i$ is
    $D/\Z/m$ with action
    $e^{\frac{2\pi i}{m}}(z)=e^{\frac{2\pi i}{m}}z$. Suppose that the local trivialization of $L$ is
    $D\times \C/\Z/m$ with the action
    $$e^{\frac{2\pi i}{m}}(z, w)=(e^{\frac{2\pi i}{m}}z, e^{\frac{2\pi id}{m}}w).$$
    When $d=0$, we call $z_i$ a Ramond marked point. When $d>0$,
    we call $z_i$ a Neveu-Schwarz marked point.}
	\vskip 0.1in
	\noindent
	{\bf Remark 5.2.3: }{\it 
    If  $L$ is an orbifold line bundle on a smooth orbifold Riemann surfae $\tilde{\Sigma}$,
     then the sheaf of local invariant sections of $L$
    is locally free of rank one, and hence dual to a unique orbifold
    line bundle $|L|$ on $\Sigma$. $|L|$ corresponds to
    the desingularization of $L$ \cite{CR1}(Prop 4.1.2). It can be constructed as follows.

    We keep the local trivialization at other places and change it at orbifold point $z_i$ by
    a $\Z/m$-equivariant map $\Psi: D-\{0\}\times \C\rightarrow D-\{0\}\times \C$ by
    $$(z, w)\rightarrow (z^m, z^{-d}w),\leqno(5.1)$$
    where $\Z/m$ acts trivially on the second $D-\{0\}\times \C$. Then, we extend $L|_{D-\{0\}\times C}$
    to a smooth holomorphic line bundle over $\Sigma$ by the second trivialization. Since $\Z/m$ acts
    trivially, hence it is a line bundle over $\Sigma$, which we denote it by
    $|L|$. Note that if $z_i$ is a Ramond marked point, $|L|=L$
    locally. When $z_i$ is a Neveu-Schwarz marked point, $|L|$
    will be different locally. In particular, it does not have a
    canonical fiber over $z_i$.}
	\vskip 0.1in
	\noindent
	{\bf Example 5.2.4: }{\it 
      An orbifold has a natural orbifold canonical bundle, defined as
    the top wedge product of its (orbifold) cotangent bundle. When
    $\tilde{\Sigma}=(\Sigma,\z, \m)$ is an orbifold Riemann surface,
    the desingularization
    $$\pi_*(K_{\tilde{\Sigma}})=K_{\Sigma}\otimes_i
    \O(-(m_i-1) z_i).$$}
	\vskip 0.1in

    Next, we study the sections. Suppose that
    $s$ is a section of $|L|$ having local representative $g(u)$. Then, $(z, z^d g(z^m))$
    is a local section of $L$.
      Therefore, we obtain
    a section $\pi(s)\in \Omega^0(L)$ which equals to $s$ outside of orbifold points under identification (5.1).
    It is clear that if $s$ is holomorphic, so it $\pi(s)$. If
    we start from an analytic section of $L$, we can reverse the above process to obtain a section of
    $|L|$. In particular, $L, |L|$ have the isomorphic space of holomorphic section. In the same way,
    there is a map $\pi: \Omega^{0,1}(|L|)\rightarrow \Omega^{0,1}(L)$, where $\Omega^{0,1}(L)$ is
    interpreted as orbifold $(0,1)$-form with value in $L$. Suppose that $g(u)d\bar{u}$ is a local
    representative of a section of $t\in \Omega^{0,1}(|L|)$. $\pi(t)$ has local representative
    $z^kg(z^m)m\bar{z}^{m-1}d\bar{z}$. Moreover, $\pi$ induces an isomorphism from $H^1(|L|)\rightarrow H^1(L)$.

    Suppose that $L^r\cong K_{log}$. $rd=lm$ for some $l$. $d<m$ implies $l<r$. Therefore, $\frac{d}{m}=\frac{l}{r}$.
    Suppose that $s\in \Omega^0(|L|)$ with local representative $g(u)$. Then, $s^r$ has local representative
    $z^{rk}g^r(z^d)=z^{ml}g^r(z^d)=u^lg^r(u)$. Hence, $s^r\in \Omega^0(K_{log} \otimes \O ((-l_i) z_i).$
    When $l_i>0$, $s^r\in \Omega^0(K).$
	\vskip 0.1in
	\noindent
	{\bf Definition 5.2.5: }{\it 
    Suppose that all the marked points are Neveu-Schwarz. We have
    $$\bar{s}^{r-1}\in \Omega^0(\bar{K}\otimes |\bar{L}|^{-1})\cong \Omega^{0,1}(|L|),$$
    Hence, $\bar{\partial}s + \bar{s}^{r-1}$ is an element of
    $\Omega^{0,1}(|L|)$. We define {\em spin equation} to be
    $$\bar{\partial}s + \bar{s}^{r-1}=0.$$}
	\vskip 0.1in

    The following Lemma is due to Witten.
    \vskip 0.1in
	\noindent
	{\bf Lemma 5.2.6: }{\it 
    Suppose that all the marked points are Neveu-Schwarz for $L$. Then, $\bar{\partial}s+\bar{s}^{r-1}=0$ iff
    $s=0$.}
	\vskip 0.1in

    {\bf Proof: }
    Fix a K\"{a}hler metric on $\Sigma$ and an Hermitian metric on $L$.
    $$\bar{\partial}s, \bar{s}^{r-1})=\int_{\Sigma} \partial
    \bar{s}\cdot \bar{s}^{r-1}=\frac{1}{r}\int_{\Sigma}\partial
    (\bar{s}^r)=0.$$
    Hence,
    $$(\bar{\partial}s+\bar{s}^{r-1},
    \bar{\partial}s+\bar{s}^{r-1})=(\bar{\partial}s,\bar{\partial}s)+
    (\bar{s}^{r-1}, \bar{s}^{r-1})=0 $$
    iff
    $$\bar{\partial}s=\bar{s}^{r-1}=0,$$
    and hence $s=0$.

    When $L$ posses a Ramond marked point $z_i$, $\bar{s}^{r-1}\in \Omega^0(\bar{K}_{log}\otimes L)
    =\Omega^{0,1}(\Sigma-\{z_i\}, L)=\Omega^{0,1}(\Sigma-\{z_i\},
    |L|)$. More generally, let $\check{\Sigma}$ be the puncture
    Riemann surface obtained by removing all the Ramond marked
    points.
    \vskip 0.1in
	\noindent
	{\bf Definition 5.2.7: }{\it 
        For $s\in \Omega^0(\check{\Sigma}, |L|)$,
    $$\bar{\partial}s+\bar{s}^{r-1}\in
    \Omega^{0,1}(\check{\Sigma}, |L|),\leqno(5.2.1)$$
    and we define spin equation to be
    $$\bar{\partial}s+\bar{s}^{r-1}=0.$$}
	\vskip 0.1in

    To have an elliptic theory, we should perform a change of
    coordinate $z\rightarrow -log z$ around puncture point and
    view $\check{\Sigma}$ as a Riemann surface with cylindric
    ends.
       It is not clear if the solution of spin equation is zero in Ramond case.

    We spend the rest of section to extend above discussion to
    nodal curve

    \vskip 0.1in
	\noindent
	{\bf Definition 5.2.8: }{\it
    A  nodal curve with $k$ marked points is a pair $(\Sigma,\z)$ of
    a connected topological space
    $\Sigma=\bigcup\pi_{\Sigma_\nu}(\Sigma_\nu)$, where $\Sigma_\nu$
    is a smooth complex curve,  and $\pi_\nu:\Sigma_\nu\rightarrow
    \Sigma$ is a continuous map, and $\z=(z_1,\cdots,z_k)$ are
    distinct $k$ points in $\Sigma$ with the following properties:
    \begin{itemize}
    \item For each $z\in\Sigma_\nu$, there is a neighborhood of it such that
    the restriction of $\pi_\nu:\Sigma_\nu\rightarrow \Sigma$ to this
    neighborhood is a homeomorphism to its image.
    \item For each $z\in\Sigma$, we have $\sum_\nu \#\pi_\nu^{-1}(z)\leq 2$.
    \item $\sum_\nu \#\pi_\nu^{-1}(z_i)=1$ for each $z_i\in\z$.
    \item The number of complex curves $\Sigma_\nu$ is finite.
    \item The set of nodal points $\{z|\sum_\nu \#\pi_\nu^{-1}(z)=2\}$ is  finite.
    \end{itemize}}
	\vskip 0.1in

    A point $z\in\Sigma_\nu$ is called {\em singular} (or a
    \emph{node}) if $\sum_\omega\#\pi_\omega^{-1}(\pi_\nu(z))=2$. A
    point $z\in\Sigma_\nu$ is said to be a {\em marked point} if
    $\pi_\nu(z)=z_i\in\z$. Each $\Sigma_\nu$ is called a {\em
    component} of $\Sigma$. Let $k_\nu$ be the number of points on
    $\Sigma_\nu$ which are either singular or marked, and $g_\nu$ be
    the genus of $\Sigma_\nu$, a nodal curve $(\Sigma,\z)$ is called
    {\em stable} if $k_\nu+2g_\nu\geq 3$ holds for each component
    $\Sigma_\nu$ of $\Sigma$.

    A map $\vartheta:\Sigma\rightarrow\Sigma^\prime$ between two nodal
    curves is called an {\em  isomorphism} if it is a homeomorphism,
    and if it can be lifted to biholomorphisms
    $\vartheta_{\nu\omega}:\Sigma_\nu\rightarrow \Sigma^\prime_\omega$
    for each component $\Sigma_\nu$ of $\Sigma$. If $\Sigma$
    ,$\Sigma^\prime$ have marked points $\z=(z_1,\cdots,z_k)$ and
    $\z^\prime=(z_1^\prime,\cdots,z_k^\prime)$, respectively, then we
    require $\vartheta(z_i) =z_i^\prime$ for each $i$. Let
    $aut(\Sigma,\z)$ be the group of automorphisms of $(\Sigma,\z)$.

    \vskip 0.1in
	\noindent
	{\bf Definition 5.2.9: }{\it A \emph{nodal
    orbicurve} is a nodal marked  curve $(\Sigma,\z)$ with an orbifold
    structure as follows:
    \begin{itemize}
    \item The set $\z_\nu$ of orbifold points of each
    component $\Sigma_\nu$ is contained in the set of marked points
    and nodal points $\z$.
    \item A neighborhood of a marked point
    is uniformized by a branched covering map
    $z\rightarrow z^{m_i}$ with $m_i\geq 1$.
    \item A neighborhood of a nodal point
    (viewed as a neighborhood of the origin of $\{x y=0\}\subset \C^2$)
    is uniformized by a branched covering map $(x,y)\rightarrow
    (x^{n_j}, y^{n_j})$, with $n_j\geq 1$, and with group action
    $e^{2 \pi i /n_j}(x,y)=(e^{2 \pi i /n_j}x,
    e^{- 2\pi i/n_j}y)$.
    \end{itemize}
    Here $m_i$ and $n_j$ are allowed to be equal to one, i.e., the
    corresponding orbifold structure is trivial there. We denote the
    corresponding nodal orbicurve by $(\Sigma,\z,\m,\n)$ where
    $\m=(m_1,\cdots,m_k)$ and $\n=(n_j)$.}
	\vskip 0.1in
	\noindent
	{\bf Definition 5.2.10: }{\it 
    An {\em  isomorphism} between two nodal orbicurves
    $\tilde{\vartheta}: (\Sigma,\z,\m,\n)\rightarrow
    (\Sigma^\prime,\z^\prime,\m^\prime,\n^\prime)$ is a collection of
    $C^\infty$ isomorphisms $\tilde{\vartheta}_{\nu\omega}$ between
    orbicurves $\Sigma_\nu$ and $\Sigma_\omega^\prime$ which induces
    an isomorphism $\vartheta:(\Sigma,\z)\rightarrow
    (\Sigma^\prime,\z^\prime)$. The {\em  group of automorphisms} of a
    nodal orbicurve $(\Sigma,\z,\m,\n)$ is denoted by
    $aut(\Sigma,\z,\m,\n)$. It is easily seen that
    $aut(\Sigma,\z,\m,\n)$ is a subgroup of $aut(\Sigma,\z)$ of
    finite index.}
	\vskip 0.1in
	\noindent
	{\bf Definition 5.2.11: }{\it  A \emph{nodal $r$-spin orbicurve} is a nodal
    orbicurve $\tilde{\Sigma}= (\Sigma, \z,\mathbf{m})$ together with
    a pair $(L,\phi)$, where $L$ is an orbifold line bundle on
    $\tilde{\Sigma}$, and $\phi: L^r\rightarrow K_{log}$ is an
    isomorphism. Suppose that the local chart of $L$ at nodal
    point is $D_0\times \C/\Z/m$, where $D_0=\{xy=0,
    |x|,|y|<\epsilon\}$ and the action is $\lambda(x, y,
    w)=(e^{\frac{2\pi i}{m}}x, e^{-\frac{2\pi i}{m}}y, e^{\frac{2\pi
    id}{m}}w)$. We call the nodal point {\em a Ramond nodal point}
    of $d=0$. Otherwise, we call it {\em a Neveu-Schwarz nodal
    point}.}
    \vskip 0.1in
	\noindent
	{\bf Remark 5.2.12: }{\it 
     Near a Ramond nodal point, the local
    invariant  sections are not necessarily locally free, but they do
    form a rank-one torsion-free sheaf (see \cite{AJ}). Hence, the desingularization
    $|L|$ is not a line bundle near a Ramond nodal point. However, we can desingularize the restriction
    of $L$ on each component and obtain a collection of line bundle over each component.
    It is clear that a NS-nodal point give rise a pair of NS-marked points on its components. Hence, the
    desingularization at this point has no canonical fiber.
    Therefore, there is no gluing condition at NS-nodal point.
    Actually, the situation is more subtle than this. For our
    purpose, it is enough to consider desingularization $|L|$ as a
    collection of line bundles over each component glued together
    at Ramond nodal points.

    Next, we define spin equation. For simplicity, we assume
    $\Sigma=\Sigma_1\vee \Sigma_2$ at nodal point $Q=\{p=q\}$.
    When $Q$ is NS,
    $$\Omega^0(|L|)=\Omega^0(|L_1|)\times \Omega^0(|L_2|), \Omega^{0,1}(|L|)=\Omega^{0,1}(|L_1|)\times \Omega^{0,1}(|L_2|)$$
    When $Q$ is Ramond,
    $$\Omega^0(|L|)=\{s_1(p)=s_2(q); s_i\Omega^0(|L_i|)\}, \Omega^{0,1}(|L|)=\{\omega_1(p)=\omega_2(q); \omega_i\in \Omega^0(\check{\Sigma}_i, |L_i|)\}$$
    In either case, there is a map
    $$\pi: \Omega^0(|L|)\rightarrow  \Omega^0(L),
    \Omega^{0,1}(|L|)\rightarrow \Omega^{0,1}(L)$$
    inducing an isomorphism on the space of holomorphic sections (forms). To see
    this, recall that $\pi(s)$ for $s\in \Omega^0(|L|)$ is
    automatically zero at NS-marked points. Hence, $\pi(s)$ can be
    glued (as zero) together at NS-nodal points even though $s$
    doesn't. For the same reason, if $\omega\in \Omega^0(|L|)$,
    $\pi(\omega)$ can be viewed as a 1-form with zero residue at
    nodal point and hence can be glued together.}
   \vskip 0.1in

    Recall that for a nodal Riemann surface $\Sigma$, its
    canonical bundle $K_{\Sigma}$ restricts to $K_{log}$ on each
    component. Therefore, spin equation is well-defined no matter
    the nodal point is Ramond or Neveu-Schwarz.

    \vskip 0.1in
	\noindent
	{\bf Lemma 5.2.13: }{\it 
    Suppose that $\Sigma$ is a nodal orbicurve and $L$ is $r$-root of $K_{log}$ such that all the marked points are
    Neveu-Schwarz.
    Then $\bar{\partial}s  +\bar{s}^{r-1}=0$ iff $s=0$.}
	\vskip 0.1in
    {\bf Proof: }
    Without the loss of generality, we assume that
    $\Sigma=\Sigma_1\vee \Sigma_2$ has only two component and no
    other marked point. Suppose that $s=(s_1, s_2)$. If the nodal
    point is NS, The lemma follows from Witten Lemma directly. Suppose
    that nodal point is Ramond,
    $$    (\bar{\partial}s, \bar{s}^{r-1})=\sum_i \int_{\Sigma} \partial
    \bar{s}_i\cdot \bar{s}_i^{r-1}=\sum_i \frac{1}{r}\int_{\Sigma}\partial
    (\bar{s}^r)=\bar{s}^r(p)-\bar{s}^r(q)=0.$$
    Hence,
    $$(\bar{\partial}s+\bar{s}^{r-1},
    \bar{\partial}s+\bar{s}^{r-1})=(\bar{\partial}s,\bar{\partial}s)+
    (\bar{s}^{r-1}, \bar{s}^{r-1})=0 $$
    iff
    $$\bar{\partial}s=\bar{s}^{r-1}=0,$$
    and hence $s=0$.

    \subsection{The moduli space of spin orbifold stable maps}

    As we show in last section, the solution of equations (5.1) can
    be described using $L$ only. In this section, we construct its
    moduli space.

	\vskip 0.1in
	\noindent
	{\bf Definition 5.3.1: }{\it  Let $(X,J)$ be an almost
    complex orbifold. An \emph{$r$-spin orbifold stable map into
    $(X,J)$} is a quadruple $(f,(\Sigma,\z,\m,\n),(L, \phi), \xi)$
    described as follows:
    \begin{enumerate}
    \item $f$ is a continuous map from a nodal
    orbicurve $(\Sigma,\z)$ into $X$ such that each
    $f_\nu=f\circ\pi_\nu$ is a pseudo-holomorphic map from
    $\Sigma_\nu$ into $X$.
    \item $\xi$ is an isomorphism class of compatible structures (see
    definition in \cite{CR2}).

    \item Let $k_\nu$ be the order of the set $\z_\nu$, namely the
    number of points on $\Sigma_\nu$ which are special (i.e. nodal or
    marked ), if $f_\nu$ is a constant map, then $2 g_\nu - 2  + k_\nu
    > 0$.

    \item At any marked or nodal point $p$ the induced homomorphism
    on the local group $\rho_{p}: G_{p}\rightarrow U(1)\times
    G_{f(p)}$ is injective.
    \end{enumerate}}

    An isomorphism between $r$-spin orbifold stable maps
    $$(f,(\Sigma,\z,\m,\n),(L, \phi), \xi) \rightarrow
    (f',(\Sigma',\z',\m',\n'),(L', \phi'), \xi')$$  
	is defined to be a
    pair $(\alpha, \beta)$ of an isomorphism $\alpha:
    (f,(\Sigma,\z,\m,\n), \xi) \rightarrow (f',(\Sigma',\z',\m',\n'),
    \xi')$ of the underlying orbifold stable maps and an isomorphism
    $\beta: f^* L' \rightarrow L$, which is compatible with $\phi$ and
    $f^*(\phi')$ in the obvious way.

    Since the number of automorphisms of $L$ fixing $\phi$ is finite,
    and since the underlying map of an $r$-spin orbifold stable map is
    a stable map, its automorphism group is finite.

    Any orbifold line bundle gives a natural representation over an orbifold Riemann surface $\Sigma$
    can be modeled locally near a non-nodal point $x_0$ as $(D\times
    \C)/(\Z/k)$, where $\Z/k$ acts on the disc  $D$ by standard
    complex multiplication by $e^{2\pi i/k}$, and on $\C$ by a
    representation $\Z/k\rightarrow U(1)$. Similarly, near a node $p$,
    $\Sigma$ can be modeled as $(D \vee D' \times \mathbf{C})/(\Z/k)$
    where $D \vee D'$ is the union of two discs $D$ and $D'$ joined at
    $0$.  $\Z/k$ acts on $D$ via multiplication by $e^{2 \pi i/k}$ and
    on $D'$ via multiplication by $e^{-2\pi i/k}$; and on $\mathbf{C}$
    by a representation $\Z/k \rightarrow U(1)$. A good orbifold map
    $f: \Sigma\rightarrow X$ induces a homomorphism $\Z/k \rightarrow
    G_{f(x_0)}$. We use $\rho_{p}$ to denote the product homomorphism
    $\Z/k\rightarrow U(1)\times G_{f(p)}$.

    If $L$ is a $r$-spin structure, $\rho_p$ has image
    in $\Z/r\times G_{f(p)}$.

    Given an $r$-spin orbifold stable map $(f,(\Sigma,\z),(L, \phi),
    \xi)$,  there is an associated homology class $f_\ast([\Sigma])$
    in $H_2(X;\Z)$ defined by
    $f_\ast([\Sigma])=\sum_\nu(f\circ\pi_\nu)_\ast[\Sigma_\nu]$. By previous argument,
     for each marked point $z$ on $\Sigma_\nu$, say
    $\pi_\nu(z) =z_i\in\z$, $\xi_\nu$ determines, by the group
    homomorphism $\rho_{z_i}$, a conjugacy class $(g_i,l_i)$ where
    $g_i\in G_{f(Z_i)}$ and $l_i \in \Z/r$. Let
    $\hat{X}$ be the disjoint union of $X_{(g, l)}$.  We have a map $ev$ sending each (equivalence class of)
    stable maps into $\hat{X}^k$ by $(f,(\Sigma,\z),\xi)\rightarrow
    ((f(z_1),(g_1,l_1)),\cdots,(f(z_k),(g_k,l_k)))$.
    $\hat{X}$ plays a
    role in r-spin orbifold theory analogous to the role played by inertial orbifold
    $\tilde{X}$ in ``usual'' (without spin structure) orbifold theory.

    Similarly, we will call a type $\hat{\x}= \prod_i X_{(g_i,l_i)}$
    \emph{Ramond} if any $l_i$ is zero, and \emph{Neveu-Schwarz}
    otherwise.

    \vskip 0.1in
    \noindent
    {\bf Remark 5.3.2: }{\it
    Despite some superficial resemblance to the moduli of orbifold
    stable maps into $X \times \B (\Z/r)$, we will show in
    Section 5.3 that the moduli of $r$-spin orbifold stable
    maps has a completely different virtual fundamental class (even of
    a different dimension) and thus different Gromov-Witten
    invariants.}
    \vskip 0.1in
    \noindent
    {\bf Definition 5.3.3: }{\it  An $r$-spin stable map
    $(f,(\Sigma,\z),(L, \phi), \xi)$ is said to be \emph{of type
    $\hat{\x}$} if $ev((f,(\Sigma,\z), \xi))\in \hat{\x}$. Given a
    homology class $A\in H_2(X;\Z)$, we let
    $\overline{\M}_{g,k,r}(X,J,A,\hat{\x})$ denote the moduli space of
    equivalence classes of $r$-spin orbifold stable maps of genus $g$,
    with $k$ marked points, homology class $A$, and type $\hat{\x}$,
    i.e., $$ \overline{\M}_{g,k,r}(X,J,A,\hat{\x})= $$
    $$\{[(f,(\Sigma,\z), (L, \phi), \xi)]|g_\Sigma=g,\#\z=k,
    f_\ast([\Sigma])=A, ev((f,(\Sigma,\z),(L,p),\xi))\in \hat{\x}\}.
    $$}
    \vskip 0.1in
    \noindent
    {\bf Remark 5.3.4: }{\it
    \begin{enumerate}
    \item When $\overline{\M}_{g,k,r}(X,J,A,\hat{\x})$ is non-empty,
    the map forgetting $(L, \phi)$ makes
    $\overline{\M}_{g,k,r}(X,J,A,\hat{\x})$ into a finite branched
    cover of the moduli space of orbifold stable maps
    $\overline{\M}_{g,k}(X,J,A,\x)$.  Here $\x$ is the image of
    $\hat{\x}$ under the map $\hat{X}^k \rightarrow \tilde{X}^k$
    forgetting the $l_i$.
    \item When $X$ is a point with no orbifold structure,
    we obtain the usual moduli space of
    $r$-spin curves \cite{J2,AJ} (Mapping to a Point).
    \item The orbifold line bundle $L$ has \emph{degree} (or first
    Chern number) $$\deg(L)=\deg(K_{log})/r =\frac{2g-2+k}{r},$$ which
    must be congruent to $ \sum \frac{l_i}{r} (\mbox{mod } \Z)$, by
    \cite[Prop 4.1.2]{CR1}.  Thus,
    $\overline{\M}_{g,k,r}(X,J,A,\mathbf{x})$ is empty unless
    $2g-2+k-\sum l_i=0 (\mbox{mod } r).$

    \end{enumerate}
    }
    \vskip 0.1in

    \subsection{r-Spin orbifold Gromov-Witten invariants}

    Following the procedure in \cite{CR2} for the moduli space of
    orbifold stable maps, we can define a natural stratification and
    topology for $\overline{\M}_{g,k,r}(X,J,A,\hat{\x})$.
    \vskip 0.1in
    \noindent
    {\bf Proposition 5.4.1: }{\it
    Suppose $X$ is either a symplectic orbifold with a symplectic form
    $\omega$ and an $\omega$-compatible almost-complex structure $J$,
    or a projective orbifold with an integrable almost-complex
    structure $J$. Then the moduli space
    $\overline{\M}_{g,k,r}(X,J,A,\hat{\x})$ is compact and
    metrizable.}

    By a slight modification of the construction in \cite{CR2}, one
    can obtain a virtual fundamental cycle for the moduli space of
    spin orbifold stable maps.
    \vskip 0.1in
    \noindent
    {\bf Theorem 5.14: }{\it If all $l_i$ are non-zero, we can construct a
    Kuranishi structure of
    $\overline{\M}_{g,k,r}(X, A, J, \hat{\x})$ in the sense of \cite{FO}. Hence, it defines a
    virtual fundamental cycle $[\overline{\M}_{g,k,r}(X, A,
    J,\hat{\x})]^{vir}$ of degree 
	$$ = 2 \left [
    c_1(A)+(n-3)(1-g)+k-\iota(\hat{\x})-\frac{1}{r}\left ( (g-1)(r-2) +
    k - \sum_i l_i \right ) \right ].$$}
    \vskip 0.1in

    For any component $\hat{\x}=(X_{(g_1,l_1)},\cdots,X_{(g_k,l_k)})$,
    there are $k$ evaluation maps $$
    e_i:\overline{\M}_{g,k,r}(X,J,A,\hat{\x})\rightarrow
    X_{(g_i,l_i)}, \hspace{4em} i=1,\cdots, k. $$ We also have a map
    $$p: \overline{\M}_{g,k,r}(X,J,A,\hat{\x})\rightarrow
    \overline{\M}_{g,k},$$ where $p$ contracts the unstable components
    of the domain to obtain a stable Riemann surface in
    $\overline{\M}_{g,k}$.

    In addition to the forgetful morphism $\overline{\M}_{g,k,r}
    \rightarrow \overline{\M}_{g,k}$, there is another important
    morphism of the moduli space of r-spin orbifold stable maps;
    namely, \emph{forgetting tails}.   When a nodal r-spin orbicurve
    $(\tilde{\Sigma},(L,\varphi))=((\Sigma, \z, \m),(L,\varphi))$ has
    a point $(z_i,m_i)$ such that the representation $\Z/m_i
    \rightarrow \Z/r \hookrightarrow U(1)$, then
    $(L,\varphi)$ induces an r-spin structure $(\pi_*L, \pi_*
    \varphi)$ on the orbicurve $\overline{\Sigma}:=(\Sigma, (z_1,
    \dots, \hat{z}_i, \dots, z_k), (m_1, \dots, \hat{m}_i, \dots,
    m_k))$ obtained by forgetting the $i$th marked point and its
    orbifold structure.  Here $\pi:\tilde{\Sigma} \rightarrow
    \overline{\Sigma}$ is the obvious ``forgetful'' map.

    Indeed, if $\hat{\Sigma}$ is the marked orbicurve
    $$\hat{\Sigma}=(\Sigma, \z, (m_1, \dots, 1, \dots, m_k))$$
    obtained from $\tilde{\Sigma}$ by forgetting the orbifold
    structure at $z_i$, but keeping the marked point $z_i$, we have a
    commuting diagram $$
    \begin{array}{rl}
    \tilde{\Sigma}  & \stackrel{P}{\rightarrow}  \hat{\Sigma}\\  ^\pi
    \downarrow &   \swarrow_q \\ \overline{\Sigma}
    \end{array}
    $$

    By Remark~\ref{rem:zeros} we have $p_* \varphi:(p_*L)^{\otimes r}
    \tilde{\rightarrow} K_{\tilde{\Sigma},\log} \otimes \O (z_i) \cong
    K_{\overline{\Sigma},\log},$ and of course $p_*L = \pi_* L$, where
    we have simply forgotten the $i$th point $z_i$.

    Thus, if the type $\hat{\x}$ corresponds to points of $\hat{X}^k$
    with $l_i=r-1$, there exists a \emph{forgetting tails} morphism
    $$\overline{\M}_{g,k,r}(X,J,A,\hat{\x}) \rightarrow
    \overline{\M}_{g,k-1,r}(X,J,A,\hat{\x}')$$ where $\hat{\x}'$ is
    the connected component of $\hat{X}^{k-1}$ obtained by mapping the
    component $\hat{\x}$ to $\hat{X}^{k-1}$, via the map forgetting
    the $i$th component of $\hat{X}^k$.

    For any set of cohomology classes $\gamma_i\in
    H^{*-2\iota(g_i)}(X_{(g_i,l_i)};\Q)\subset H^*_{orb}(X;\Q)$,
    $i=1,\cdots,k$ and $\sigma \in H^*(\overline{\M}_{g,k}, \Q)$, the
    $r$-spin orbifold Gromov-Witten invariant is defined as the
    pairing $$ \Psi^{X,J}_{(g,k,r,A,\hat{\x})}(\sigma;
    \tau_{a_1}(\gamma_1), \cdots, \tau_{a_k}(\gamma_k))=p^*\sigma
    (\prod_{i=1}^k \psi_{i}^{a_i}) e^*_i(\sigma_i) \cap
    [\overline{\M}_{g,k,r}(X,J,A,\hat{\x})]^{vir}. $$ where $\psi_i$
    is the first Chern class of the line bundle generated by the
    cotangent space of the $i$-th marked point. When $a_i<0$, we
    define it to be zero. The same argument as in the case of ordinary
    GW-invariants yields

    \vskip 0.1in
    \noindent
    {\bf Proposition 5.4.2: }{\it
    \begin{enumerate}
    \item $ \Psi^{X,J}_{(g,k,r,A,\hat{\x})}(\sigma; \tau_{a_1}(\gamma_1), \cdots,
    \tau_{a_k}(\gamma_k))=0$ unless $\deg \sigma+ \sum_i
    (deg_{orb}(\gamma_i)+a_i)=2c_1(A)+2(n-3)(1-g)+2k-2(g-1)(1-\frac{2}{r})-\frac{2k}{r}+
    2 \sum l_i/r$ , where $deg_{orb}(\gamma_i)$ is the orbifold degree
    of $\gamma_i$ obtained after degree shifting.
    \item $\Psi^{X,J}_{(g,k,r,A,\hat{\x})}(\sigma; \tau_{a_1}(\gamma_1), \cdots,
    \tau_{a_k}(\gamma_k))$ is independent of the choice of $J$ and
    hence is an invariant of symplectic structures.
    \end{enumerate}}
    \vskip 0.1in

    Spin orbifold Gromov-Witten invariants are expected to satisfy
    a set of interesting axioms. Some of them are similar to the
    axiom of ordinary Gromov-Witten invariants. But some are new
    addition with the introduction of spin curves. We are in the
    process to verify these axioms.

\section{Conclusion Remark}
    During the last year and half, this new subject is experiencing an
    explosion. This author often hears that a new idea would run through the
    physics like wild fire. We mathematician are used to be on the
    slower pace. However, it is fair to say that in such a short
    time, the basic idea of stringy orbifold is spreading like a wild
    fire in mathematics. This speed of new results is even beyond this
    author's expectation. The excitement  reminds him a lot the
    quantum cohomology era during 94-96. The diversity of papers in this volume is
    the evidence of such an explosion.

    Because of the limited time and space, this author did not
    include several topics which he initially intends to.
    Recently, there is an algebraic construction of the orbifold
    quantum cohomology over Deligne-Mumford stack including an
    interesting integral version of the orbifold cohomology.
    Furthermore, the construction of the orbifold cohomology of
    \cite{CR1} has been adapted to construct an orbifold Chow rings
    \cite{AV1}, \cite{F}. This exciting development is surveyed by
    their authors in this volume.

    There are some very interesting constructions of the orbifold
    elliptic genus \cite{BL}, \cite{DLM}. Physically, it belongs to a
    different type of string theory (heterotic string theory).
    Mathematically, they are also quite different from other
    orbifold theory such as orbifold cohomology. For examples, the
    twisted sectors in the orbifold cohomology is based on the commuting
    pair while the twisted sector in the elliptic genus is based on
    commuting triple.

    There is very recent work of Borisov-Mavlyutov \cite[BM] who are considering
    the mirror symmetry from the vertax algebra point of view. In the
    process, they seems to suggest some conjectural solution for the 
    orbifold cohomology product for certain toric varieties. It
    would be an interesting problem to check their answer.

    There are also very recent papers by Lupercio-Poddar \cite{LP}
    and Yasuda \cite{Y}, where a weak form of the orbifold string theory
    conjecture for the equivalence of orbifold the Hodge number and
    Hodge number of the crepant resolution was solved in its complete
    generality. This is an exciting development.

    We are living in a good time. This author expects more excitements to
    come!


\begin{thebibliography}{L}
    \bibitem[AJ]{AJ} D.~Abramovich, T.~Jarvis, {\em Moduli of twisted spin curves},
     \texttt{math.AG/0104154}.

    \bibitem[AV]{AV}
    D.~Abramovich, A.~Vistoli, {\em Compactifying the space of stable
    maps}, \texttt{math.AG/9908167}.

    \bibitem [AR]{AR} A. Adem and Y. Ruan, {\em Twisted orbifold
    K-theory}, math.AT/0107168
    \bibitem [AV1]{AV1} D. Abramovich and A. Vistoli, {\em Lectures
    on Madison conference},
    \bibitem [A]{A} P. Aspinwall, {\em Enhanced gauge symmetries and K3-surface}, Phys Lett B357(1995)329-334
    \bibitem [B]{B} V. Batyrev, {\em Birational Calabi-Yau n-folds have the equal betti numbers}, math.AG/9710020
    \bibitem [B1]{B1} V. Batyrev, {\em Mirror symmetry and toric
    geometry,} Doc. Math. Vol II 239-249
    \bibitem[Bea]{Bea}  A. Beauville,
{\it Sur la cohomologie de certains espaces de modules de fibr\'es
vectoriels},  Geometry and Analysis (Bombay, 1992), 37-40, Tata
Inst. Fund. Res., Bambay, 1995.
    \bibitem [BL]{BL} L. Borisov and A. Libgober, {\em Elliptic
    genera of singular varieties}, math.AG/0007108
    \bibitem [BM]{BM} L. Borisov and A. Mavlyutov, {\em String cohomology of Calabi-Yau hypersurfaces via Mirror
    Symmetry}, math.AG/0109096
    \bibitem [C]{C} W. Chen, {\em A homotopy theory of
    orbispaces}, math.AT/0102020
    \bibitem[DHVW]{DHVW} L. Dixon, J. Harvey, C. Vafa and E.
    Witten, {\em Strings on orbifolds, I, II}, Nucl. Phys.
    B261(1985), 678, B274(1986), 285.
    \bibitem [DLM]{DLM} C. Dong, K. Liu and X. Ma, {\em On orbifold elliptic
    genus}, math.DG/0109005
    \bibitem [CM]{CM} M. Crainic and I. Moerdijk, {\em A homology
    theory for \'etale groupoids}, J. Reine. Angew. Math.
    521(2000), 25-46.
    \bibitem [CR1]{CR1} W. Chen and Y. Ruan, {\em A new cohomology theory
    for orbifold}, math.AG/0004129
    \bibitem [CR2]{CR2} W. Chen and Y. Ruan, {\em Orbifold Gromov-Witten theory}, math.AG/0103156
    \bibitem [D]{D} Dijkgraaf, R. \emph{Discrete Torsion and Symmetric
Products}, preprint hep-th/9912101.
    \bibitem [DL]{DL} J. Denef and F. Loeser, {\em Motivic integration, quotient singularities
    and the McKay correspondence,} math.AG/9903187
    \bibitem[EGL]{EGL} G. Ellingsrud, L. G\"ottsche,  M. Lehn,  {\em On the
cobordism class of the Hilbert scheme of a surface}, J. Algebraic
Geom. {\bf 10} (2001) 81--100.

\bibitem[ES1]{ES1}  G. Ellingsrud,   S.A. Str\o mme,
{\em On the homology of the Hilbert scheme of points in the
plane}, Invent. Math.  {\bf 87} (1987) 343--352.

\bibitem[ES2]{ES2} G. Ellingsrud,  S. A. Str\o mme,
{\em Towards the Chow ring of the Hilbert scheme of $\P^2$}, J.
reine angew. Math. {\bf 441} (1993) 33--44.

\bibitem [F]{F} B. Fantechi, in preparation.
\bibitem[FG1]{FG1}  B. Fantechi, L. G\"ottsche,
{\em The cohomology ring of the Hilbert schemes of $3$ points on a
smooth projective variety}, J. reine angew. Math.  {\bf 439}
(1993) 147--158.
        \bibitem [FG2]{FG2} B. Fantechi and L. G\"{o}ttsche, {\em
    Orbifold cohomology for global quotients}, math.AG/0104207
    \bibitem [FO]{FO} K. Fukaya and K. Ono, {\em Arnold conjecture and
    Gromov-Witten invariant}, Topology {\textbf 38}(1999), no.5,
    933-1048.
    \bibitem [FR]{FR} D. Freed, {\em The Verlinde algebra is twisted
    equivariant  K-theory}, math.RT/0101038
    \bibitem[FW]{FW} I. Frenkel and W. Wang,
{\em Virasoro algebra and wreath product convolution}, J.~Alg.
{\bf 242} (2001) 656--671.
    \bibitem [G]{G} L. G\"{o}ttsche, {\em The Betti numbers of the Hilbert scheme of points
               on a smooth projective surface}. Math. Ann. 286 (1990), no. 1-3, 193--207.
    \bibitem[Gro]{Gro} I.~Grojnowski,
{\em Instantons and affine algebras I: the Hilbert scheme and
vertex operators}, Math. Res. Lett. {\bf 3} (1996) 275--291.
    \bibitem[H]{H} Haeberly, J.-P. \emph{For $G=S^1$
there is no $G$-Chern Character}, Contemporary Mathematics
\textbf{36} (1985), pp. ~113--118.
    \bibitem[J1]{J1} T.~J.~Jarvis, {\em Geometry of the moduli of higher spin
    curves,} Internat. J. of Math. \textbf{11} (2000), no. 5.
    \texttt{math.AG/9809138}.

    \bibitem[J2]{J2} T.~J.~Jarvis, {\em Torsion-free sheaves and moduli of
    generalized spin curves,} Compositio Math. \textbf{110} (1998) no. 3,
    291-333.

    \bibitem[JKV]{JKV} T.~Jarvis, T.~Kimura, A.~Vaintrob, {\em Moduli spaces of higher
    spin curves and integrable hierarchies}, Compositio Math. \textbf{126} (2001), no. 2,
    157--212. \texttt{math.AG/9905034}.

    \bibitem[JKV2]{JKV2} T.~Jarvis, T.~Kimura, A.~Vaintrob, {\em Gravitational descendants and the
    moduli space of higher spin curves.}  Advances in algebraic geometry
    motivated by physics (Lowell, MA, 2000), Contemp. Math., \textbf{276},
    Amer. Math. Soc., Providence, RI, 2001, 167--177. \texttt{math.AG/9905034}.

    \bibitem [JR]{JR} T. Jarvis and Y. Ruan, {\em Spin orbifold
    quantum cohomology}, in preparation.

    \bibitem[K]{K}
    M.~Kontsevich, {\em Intersection theory on the moduli space of
    curves and the matrix Airy function,} Commun. Math. Phys.
    \textbf{147} (1992), 1--23.

    \bibitem[Leh]{Leh} M. Lehn,
{\em Chern classes of tautological sheaves on Hilbert schemes of
points on surfaces}, Invent. Math. {\bf 136} (1999) 157--207.


\bibitem[LS1]{LS1} M. Lehn and C. Sorger, {\em Symmetric groups and
the cup product on the cohomology of Hilbert schemes}, Duke Math.
J. (to appear), math.AG/0009131.


    \bibitem [LS2]{LS2} M. Lehn and C. Sorger, {\em The cup product of the Hilbert scheme for
    $K3$-surfaces}, math.AG/0012166
    \bibitem [LQ]{LQ} W. Li and Z. Qin, {\em Toward the
    Gromov-Witten invariants of the Hilbert schemes of points in
    algebraic surfaces}, preprint
    \bibitem[LQW1]{LQW1} W.-P. Li, Z. Qin and W. Wang, {\em Vertex algebras and the
cohomology ring structure of Hilbert schemes of points on
surfaces}, Math. Ann. (to appear), math.AG/0009132.


\bibitem[LQW2]{LQW2} W.-P. Li, Z. Qin and W. Wang, {\em Generators
for the cohomology ring of Hilbert schemes of points on surfaces},
Internat. Math. Res. Notices (to appear), math.AG/0009167.


\bibitem[LQW3]{LQW3} W.-P. Li, Z. Qin and W. Wang, {\em Universality and
stability of cohomology rings of Hilbert schemes of points on
surfaces}, Preprint, math.AG/0107139.

\bibitem[LQW4]{LQW4} W.-P. Li, Z. Qin and W. Wang, {\em Hilbert schemes and $W$
algebras}, Preprint.


    \bibitem [LR]{LR} A. Li and Y. Ruan, {\em Symplectic surgery
    and Gromov-Witten invariants of Calabi-Yau 3-folds I}, to
    apear in Invent. Math.
    \bibitem [LP]{LP} E. Lupercio and M. Poddar, in preparation.
    \bibitem [LU1]{LU1} E. Lupercio and B. Uribe, {\em Gerbes over
    orbifolds and twisted K-theory}, math.AT/0105039
    \bibitem [LU2]{LU2} E. Lupercio and B. Uribe, {\em Loop
    groupoids, gerbes and twisted sectors on orbifolds,}
    math.AT/0110207
    \bibitem[Mar]{Mar} E. Markman,
{\em  Generators of the cohomology ring of moduli spaces of
sheaves on symplectic surface}, Preprint, math.AG/0009109.
    \bibitem [MP1]{MP1} I. Moerdijk and D. Pronk, {\em Orbifolds,
    sheaves and groupoids}, K-theory 12 (1997), no 1, 3-21
    \bibitem [MP2]{MP2} I. Moerdijk and D. Pronk, {\em Simplicial
    cohomology of orbifolds}, Indag. Math. (N.S.) 10 (1999), no.2,
    269-293.
    \bibitem [N]{N} H. Nakajima, {\em Lectures on Hilbert schemes of points on surfaces}.
    University Lecture Series, 18.
    \bibitem [NW]{NW}W. Nahm and K. Wendland, {\em Mirror Symmetry on Kummer Type surfaces}, private communication.
    \bibitem [QW]{QW} Z. Qin and W. Wang, {\em Hilbert schemes and symmetric products: 
	a dictionary}, appear on Proceedings of Workshop of Mathematical Aspect of Orbifold String
	Theory.
	\bibitem [R1]{R1} Y. Ruan, {\em Surgery, quantum cohomology and birational
    geometry}, AMS Translations, Ser 2. Vol 196, 183-199.
    \bibitem [R2]{R2} Y. Ruan, {\em Stringy geometry and topology
    of orbifolds}, math.AG/0011149
    \bibitem [R3]{R3} Y. Ruan, {\em Discrete torsion and twisted
    orbifold cohomology}, math. AG/0005299
    \bibitem [R4]{R4} Y. Ruan, {\em Cohomology ring of crepant
    resolutions of orbifolds}, math.AG/018195
    \bibitem [S]{S} I. Satake, {\em The Gauss-Bonnet theorem for
    V-manifolds}, J. Math. Soc. Japan 9 (1957), 464-492.
    \bibitem [U]{U} B. Uribe, {\em Orbifold cohomology of symmetry
    product}, math.AT/0109125
    \bibitem [U1]{U1} B. Uribe, Wisconsin Ph.D thesis.
    \bibitem[Vas]{Vas} E. Vasserot,
{\em  Sur l'anneau de cohomologie du sch\'ema de Hilbert de
$\C^2$}, C. R. Acad. Sci. Paris Ser. I Math. {\bf 332} (2001)
7--12.
    \bibitem [VW]{VW} C. Vafa and E. Witten, {\em A strong coupling
    test of S-duality}, Nucl Phys B 431 (1994) 3-77
    \bibitem [W]{W} W. Wang, {\em Equivariant K-theory,
    generalized symmetry product and twisted Heisenberg algebra},
    math.QA/0104168

\bibitem [WA]{WA} C-L. Wang, {\em On the topology of birational
minimal models}, J. Diff. Geom. 50 (1998), 129-146.
    \bibitem [WE]{WE} Katrin Wendland, {\em Consistency of orbifold conformal field theories on $K3$.}
                       hep-th/0010281
    \bibitem [WI]{WI} P. Wilson, {\em Flops, Type III contractions and Gromov-Witten invariants on Calabi-Yau
    threefolds}, math.AG/9707008
    \bibitem[WIT1]{WIT1}
    E.~Witten, {\em Algebraic geometry associated with matrix models
    of two dimensional gravity}, Topological methods in modern
    mathematics (Stony Brook, NY, 1991), Publish or Perish, Houston,
    TX 1993, 235-269.
    \bibitem [WIT2]{WIT2}. E. Witten, {\em D-brane and K-theory},
    hep-th/9810188

     \bibitem [Y]{Y} T. Yasuda, {\em Twisted jet, motivic measure and orbifold
    cohomology}, math.AG/0110228
    \bibitem [Z]{Z} Wuanchuan Zhang, Ph.D Thesis, University of
    Wisconsin.
    \end{thebibliography}
\end{document}